\documentclass[11pt]{amsart}
\usepackage{amssymb}

\DeclareMathOperator{\dom}{Dom}
\DeclareMathOperator{\supp}{Supp}

\newcommand{\dbar}{\ensuremath{\overline\partial}}
\newcommand{\dbarstar}{\ensuremath{\overline\partial^*}}
\newcommand{\C}{\ensuremath{\mathbb{C}}}
\newcommand{\R}{\ensuremath{\mathbb{R}}}

\makeatletter
\newcommand{\sumprime}{\if@display\sideset{}{'}\sum%
            \else\sum'\fi}
\makeatother

\begin{document}

\numberwithin{equation}{section}

% define theorem environments
\newtheorem{theorem}{Theorem}[section]
\newtheorem{proposition}[theorem]{Proposition}
\newtheorem{conjecture}[theorem]{Conjecture}
\def\theconjecture{\unskip}
\newtheorem{corollary}[theorem]{Corollary}
\newtheorem{lemma}[theorem]{Lemma}
\newtheorem{observation}[theorem]{Observation}
\theoremstyle{definition}
\newtheorem{definition}{Definition}
\numberwithin{definition}{section}
\newtheorem{remark}{Remark}
\def\theremark{\unskip}
\newtheorem{question}{Question}
\def\thequestion{\unskip}
\newtheorem{example}{Example}
\def\theexample{\unskip}
\newtheorem{problem}{Problem}

\def\vvv{\ensuremath{\mid\!\mid\!\mid}}
\def\intprod{\mathbin{\lr54}}
\def\reals{{\mathbb R}}
\def\integers{{\mathbb Z}}
\def\N{{\mathbb N}}
\def\complex{{\mathbb C}\/}
\def\distance{\operatorname{distance}\,}
\def\spec{\operatorname{spec}\,}
\def\interior{\operatorname{int}\,}
\def\trace{\operatorname{tr}\,}
\def\cl{\operatorname{cl}\,}
\def\essspec{\operatorname{esspec}\,}
\def\range{\operatorname{\mathcal R}\,}
\def\kernel{\operatorname{\mathcal N}\,}
\def\linearspan{\operatorname{span}\,}
\def\lip{\operatorname{Lip}\,}
\def\sgn{\operatorname{sgn}\,}
\def\Z{ {\mathbb Z} }
\def\e{\varepsilon}
\def\p{\partial}
\def\rp{{ ^{-1} }}
\def\Re{\operatorname{Re\,} }
\def\Im{\operatorname{Im\,} }
\def\dbarb{\bar\partial_b}
\def\eps{\varepsilon}

\def\Hs{{\mathcal H}}
\def\E{{\mathcal E}}
\def\scriptu{{\mathcal U}}
\def\scriptr{{\mathcal R}}
\def\scripta{{\mathcal A}}
\def\scriptc{{\mathcal C}}
\def\scriptd{{\mathcal D}}
\def\scripti{{\mathcal I}}
\def\scriptk{{\mathcal K}}
\def\scripth{{\mathcal H}}
\def\scriptm{{\mathcal M}}
\def\scriptn{{\mathcal N}}
\def\scripte{{\mathcal E}}
\def\scriptt{{\mathcal T}}
\def\scriptr{{\mathcal R}}
\def\scripts{{\mathcal S}}
\def\scriptb{{\mathcal B}}
\def\scriptf{{\mathcal F}}
\def\scriptg{{\mathcal G}}
\def\scriptl{{\mathcal L}}
\def\scripto{{\mathfrak o}}
\def\scriptv{{\mathcal V}}
\def\frakg{{\mathfrak g}}
\def\frakG{{\mathfrak G}}

\def\ov{\overline}
%\date {November, 2004.  %Print \today}

\author{Siqi Fu}
\thanks
{The author was supported in part by NSF grants DMS
0070697/0500909 and by an AMS Centennial Research Fellowship.}
\address{Department of Mathematical Sciences,
Rutgers University-Camden, Camden, NJ 08102}
\email{sfu@camden.rutgers.edu}
\title[]  % runningheader
{Hearing  the type of a domain in $\C^2$ \\
with the $\dbar$-Neumann Laplacian} \maketitle

%\begin{center}
%\fbox{Draft } \end{center}

\tableofcontents

\

\section{Introduction}\label{intro}

Motivated by Mark Kac's famous question ``Can one hear the shape
of a drum?" (see \cite{Kac66, GordonWebbWolpert92}), we study the
interplays between the geometry of a bounded domain in $\C^n$ and
the spectrum of the $\dbar$-Neumann Laplacian. Since the work of
Kohn~\cite{Kohn63, Kohn64}, it has been discovered that regularity
of the $\dbar$-Neumann Laplacian is intimately connected to the
boundary geometry. (See, for example, the surveys
\cite{BoasStraube99, Christ99, DangeloKohn99, FuStraube01} and the
books \cite{FollandKohn72, Hormander90, Krantz92, ChenShaw99,
Ohsawa02}.) It is then natural to expect that one can ``hear" more
about the geometry of a bounded domain in $\C^n$ with the
$\dbar$-Neumann Laplacian than with the usual Dirichlet/Neumann
Laplacians.

For bounded domains in $\C^n$, it follows from H\"{o}rmander's
$L^2$-estimates of the $\dbar$-operator \cite{Hormander65}  that
pseudoconvexity implies positivity of the spectrum of the
$\dbar$-Neumann Laplacian on all $(0, q)$-forms, $1\le q\le n-1$.
The converse is also true (under the assumption that the interior of
the closure of the domain is the domain itself).  This is a
consequence of the sheaf cohomology theory dated back to Oka and H.
Cartan--it is proved, implicitly, in \cite{Serre53, Laufer66}. (See
\cite{Fu05} for a discussion and proofs of this and other facts
without the sheaf cohomology theory.) Therefore, in Kac's language,
we can ``hear" pseudoconvexity via the $\dbar$-Neumann Laplacian.

Regularity and spectral theories of the $\dbar$-Neumann Laplacian
closely intertwine. For example, on the one hand, by a classical
theorem of Hilbert in general operator theory, compactness of the
$\dbar$-Neumann operator is equivalent to emptiness of the essential
spectrum of the $\dbar$-Neumann Laplacian.  On the other hand, by a
result of Kohn and Nirenberg~\cite{KohnNirenberg65}, compactness of
the $\dbar$-Neumann operator implies exact global regularity of the
$\dbar$-Neumann Laplacian on $L^2$-Sobolev spaces. It was shown in
\cite{FuStraube98} that for a bounded convex domain in $\C^n$, the
essential spectrum of the $\dbar$-Neumann Laplacian on $(0, q)$-form
is empty if and only if the boundary contains no $q$-dimensional
complex varieties. However, such characterization does not hold even
for pseudoconvex Hartogs domains in $\C^2$ (see \cite{Matheos97},
also \cite{FuStraube01}). Recently, it was proved in
\cite{ChristFu05} that for smooth bounded pseudoconvex Hartogs
domains in $\C^2$, emptiness of the essential spectrum of
$\dbar$-Neumann Laplacian on $(0, 1)$-forms implies that the
boundary contains no pluripotentials (more precisely, it satisfies
property ($P$) in the sense of Catlin \cite{Catlin84b} or
equivalently is $B$-regular in the sense of Sibony \cite{Sibony87}).
This, together with an earlier result of Catlin \cite{Catlin84b}
(cf. \cite{Straube97}), shows that one can determine whether the
boundary of a Hartgos domain in $\C^2$ contains pluripotentials via
the spectrum of the $\dbar$-Neumann Laplacian.

In this paper, we continue our study of the spectral theory of the
$\dbar$-Neumann Laplacian.  Our main result can be stated as
follows.

\begin{theorem}\label{maintheorem}
Let $\Omega$ be a smooth bounded pseudoconvex domain in $\C^2$. Let
$\scriptn(\lambda)$ be the number of eigenvalues of the
$\dbar$-Neumann Laplacian that are less than or equal to $\lambda$.
Then $b\Omega$ is of finite type if and only if $\scriptn(\lambda)$
has at most polynomial growth.
\end{theorem}

Recall that a smooth domain is strictly pseudoconvex if for each
boundary point, there is a local change of holomorphic coordinates
which makes the boundary near this point strictly convex (in the
sense that the real Hessian of the defining function is positive
definite on the real tangent space). It is pseudoconvex if it can be
exhausted by strictly pseudoconvex ones. (We refer the reader to the
textbook \cite{Krantz01} for a treatise of these concepts.) The type
of a smooth boundary $b\Omega$ (in the sense of Kohn \cite{Kohn72})
is the maximal order of contact of a (regular) complex variety with
$b\Omega$. (See \cite{Dangelo82, Catlin84a, Dangelo93} for more
information on this and other notions of finite type.)

Theorem~\ref{maintheorem} consists of two parts. More precisely, for
the sufficiency, we establish the following result.

\begin{theorem}\label{maintheorem2a}
Let $\Omega\subset\C^2$ be smooth bounded pseudoconvex domain of
finite type $2m$. Then $\scriptn(\lambda)\lesssim \lambda^{m+1}$.
\end{theorem}

The Weyl type asymptotic formula for $\scriptn(\lambda)$ for
strictly pseudoconvex domains in $\C^n$ was established in
\cite{Metivier81} by Metivier via an analysis of the spectral kernel
of the $\dbar$-Neumann Laplacian. The heat kernel of the
$\dbar$-Neumann Laplacian on strictly pseudoconvex domains, as well
as that of the Kohn Laplacian on the boundary, were studied in a
series of papers \cite{Stanton84, BealsGreinerStanton84,
StantonTartakoff84, BealsStanton87, BealsStanton88} by various
authors. Metivier's formula was recovered as a consequence.
Recently, the heat kernel of the Kohn Laplacian on finite type
boundaries in $\C^2$ was studied by Nagel and Stein
\cite{NagelStein01}, from which one could also deduce a result
similar to Theorem~\ref{maintheorem2a} for the Kohn Laplacian on the
boundary.

We follow Metivier's approach in proving
Theorem~\ref{maintheorem2a} by studying the spectral kernel. We
are also motivated by the work on the Bergman kernel by Nagel et
al \cite{NRSW89} and by Catlin \cite{Catlin89} and McNeal
\cite{Mcneal89} as well as related work of Christ \cite{Christ88}
and Fefferman and Kohn \cite{FeffermanKohn88}. There is one
important distinction between the spectral kernel and the Bergman
kernel: While the Bergman kernel transforms well under
biholomorphic mappings, the spectral kernel does not. It does not
transform well even under non-isotropic dilations. To overcome
this difficulty, instead of (locally) rescaling the domain to unit
scale and studying the $\dbar$-Neumann Laplacian on the rescaled
domain as in the Bergman kernel case, we rescale both the domain
and the $\dbar$-Neumann Laplacian as in \cite{Metivier81} (see
Section~\ref{rescaled}). In doing so, we are led to study
non-isotropic bidiscs that have larger radii in the complex normal
direction. Roughly speaking, at a boundary point of type $2m$, the
quotient of the radii in the complex tangential and normal
directions for the bidiscs used here is $\tau:\tau^{m}$ while in
the Bergman kernel case it is $\tau:\tau^{2m}$ ($\tau>0$ is
small). To establish desirable properties, such as doubling and
engulfing properties, for these non-isotropic bidiscs, we employ
both pseudoconvexity and the finite type condition (see
Section~\ref{polydisc}). Note that only the finite type condition
was used in establishing these properties for the smaller bidiscs
used in the Bergman kernel case. Here in our analysis of these
bidiscs, we make essential use of an observation by Forn{\ae}ss
and Sibony~\cite{FornaessSibony89}. Also crucial to our analysis
is a uniform Kohn type estimate on the rescaled $\dbar$-Neumann
Laplacian. Here the sesquilinear form that defines the rescaled
$\dbar$-Neumann Laplacian controls not only the tangential Sobolev
norm of some positive order $\eps$ of any $(0, 1)$-form $u$
compactly supported on the larger bidiscs but also the square of
the quotient of the radii times the tangential Sobolev norm of
order $-1+\eps$ of the bar derivative of $u$ in the complex normal
direction (see Lemma~\ref{Q-estimate} below). By carefully
flattening the boundary, we then reduce the problem to estimating
eigenvalues of auxiliary operators on the half-space (see
Section~\ref{aux}). The final steps of the proof of
Theorem~\ref{maintheorem2a} is given in Section~\ref{lower}.

For the necessity, we prove the following slightly more general
result.

\begin{theorem}\label{maintheorem3}
Let $\Omega$ be a smooth bounded pseudoconvex domain in $\C^n$.
Let $\scriptn_q(\lambda)$ be the number of eigenvalues of the
$\dbar$-Neumann Laplacian on $(0, q)$-forms that are less than or
equal to $\lambda$.  If $\scriptn_q(\lambda)$ has at most
polynomial growth for some $q$, $1\le q\le n-1$, then  $b\Omega$
is of finite $D_{n-1}$-type.
\end{theorem}

Recall that the $D_{n-1}$-type of $b\Omega$ is the maximal order of
contact of $(n-1)$-dimensional (regular) complex varieties with
$b\Omega$.  It was observed by D'Angelo \cite{Dangelo87} that the
$D_{n-1}$-type is identical to the second entry in Catlin's
multitype. An ingredient in the proof of Theorem~\ref{maintheorem3}
is a wavelet construction of Lemari\'{e} and Meyer
\cite{LemarieMeyer86} (see Section~\ref{proof1}).  A result similar
to Theorem~\ref{maintheorem3} for the Kohn Laplacian on the
boundaries in $\C^2$ is also known to M.~Christ~\cite{Christ}.

Throughout the paper, we use $C$ to denote positive constants
which may be different in different appearances. For the reader's
convenience, we make an effort to have our presentation
self-contained.

\bigskip

\noindent{\bf Acknowledgment:} The author is indebted to Professors
M.~Christ, H.~Jacobowitz, J.~J.~Kohn,  T.~Ohsawa, N.~Stanton, and
Y.-T.~Siu for stimulating conversations and kind encouragement.

\section{Preliminaries}\label{prelim}

We first recall several relevant facts from the general operator
theory. Let $H_1$, $H_2$ be complex Hilbert spaces. Let $T$ be a
compact operator from $H_1$ into $H_2$. Then it follows from the
min-max principle that the singular values of $T$ ({\it i.e.}, the
non-zero eigenvalues of $|T|=(T^*T)^{1/2}$) is given by
\begin{equation}\label{minmax-0}
\lambda_j(T)=\inf_{g_1, \ldots, g_{j-1}\in H_1}\sup\{\|Tf\| \mid
f\perp g_1, \ldots, g_{j-1}; \|f\|=1\}, \quad j=1, 2, \ldots,
\end{equation}
where the singular values are arranged in decreasing order and
repeated according to multiplicity ({\it e.g.},
\cite{Weidmann80}).  (Throughout this paper, we will use
$\lambda_j(T)$ to denote the eigenvalues/singular values of a
compact operator $T$ arranged in this order as well as the
eigenvalues of an unbounded operator $T$ with compact resolvent
arranged in the reverse order.)  It follows that for compact
operators $T_1, T_2\colon H_1\to H_2$ and $T_3\colon H_2\to H_3$,

\begin{align}
\lambda_{j+k+1}(T_1+T_2)&\le \lambda_{j+1}(T_1)+\lambda_{j+1}(T_2)
\label{minmax-1}
\\
\intertext{and} \lambda_{j+k+1}(T_3\circ T_1) &\le
\lambda_{j+1}(T_1) \lambda_{k+1}(T_3). \label{minmax-2}
\end{align}

Let $Q$ be a non-negative, densely defined, and closed
sesquilinear form on a complex Hilbert space $H$. Then $Q$
uniquely determines a non-negative, densely defined, and
self-adjoint (unbounded) operator $S$ such that
$\dom(S^{1/2})=\dom(Q)$ and
\[
Q(u, v)=(Su,\; v)=(S^{1/2}u, \; S^{1/2}v)
\]
 for all $u\in \dom(S)$ and $v\in \dom(Q)$. (See, for example,
\cite{Davies95} for related material.) For any subspace
$L\subset\dom(Q)$, let $\lambda(L)=\sup\{Q(u, u) \mid  u\in L,
\|u\|=1\}$.  For any positive integer $j$, let
\begin{equation}\label{minmax}
\lambda_{j} (Q)=\inf\{\lambda(L) \mid L\subset \dom(Q),
\dim(L)=j\}.
\end{equation}
It follows that the associated operator $S$ has compact resolvent
if and only if $\lambda_{j}(Q)\to\infty$ as $j\to\infty$. In this
case, $\lambda_{j}(Q)$ equals $\lambda_j(S)$,  the $j^{\text
{th}}$ eigenvalue of $S$.

We now recall the setup for the $\dbar$-Neumann Laplacian ({\it
e.g.}, \cite{FollandKohn72, ChenShaw99}). Let $\Omega$ be a
bounded domain in $\C^n$. For $1\le q\le n$, let $L^2_{(0,
q)}(\Omega)$ denote the space of $(0, q)$-forms with square
integrable coefficients and with the standard Euclidean inner
product whose norm is given by
\[
\|\sumprime a_{J}d\bar z_J\|^2= \sumprime \int_\Omega |a_{J}|^2
dV(z),
\]
where the prime indicates the summation over strictly increasing
$q$-tuples $J$. (We consider $a_J$ to be defined on all $q$-tuples,
antisymmetric with respect to $J$.) For $0\le q\le n-1$, let
$\dbar_q\colon L^2_{(0, q)}(\Omega)\to L^2_{(0, q+1)}(\Omega)$ be
the $\dbar$-operator defined in the sense of distribution. This is a
closed and densely defined operator. Let $\dbarstar_q$ be its
adjoint. For $1\le q\le n-1$, let
\[
Q_q(u, v)=(\dbar_q u, \dbar_q v)+(\dbarstar_{q-1} u,
\dbarstar_{q-1} v)
\]
be the sesquilinear form on $L^2_{(0, q)}(\Omega)$ with
$\dom(Q_q)=\dom(\dbar_q)\cap\dom(\dbarstar_{q-1})$. It is evident
that $Q_q$ is non-negative, densely defined, and closed. The
operator associated with $Q_q$ is the $\dbar$-Neumann Laplacian
$\square_q$ on $L^2_{(0, q)}(\Omega)$. The following lemma is a
simple consequence of the min-max principle \eqref{minmax}.

\begin{lemma}\label{spectral}
Suppose $\lambda_{j}(Q)\ge C_1 j^\eps$ for some constants $C_1>0$
and $\eps>0$. If $u_k\in\dom(Q)$, $1\le k\le j$, satisfy
\[
\|\sum_{k=1}^j c_k u_k\|^2\ge C_2 \sum_{k=1}^j |c_k|^2
\]
for some constant $C_2>0$ and for all $(c_1, \ldots, c_j)\in
\C^j$, then
\[
\max_{1\le k\le j} Q(u_k, u_k)\ge C_1 C_2 j^\eps/(1+\eps).
\]
\end{lemma}

\begin{proof} Let $\tilde\lambda_k$, $1\le k\le j$, be the eigenvalues of
the Hermitian matrix $M=\left(Q(u_k, u_l)\right)_{1\le k, l\le
j}$. Then by the min-max principle,
\[
\tilde\lambda_k=\inf\{\tilde\lambda (\tilde L); \ \tilde L\subset
\C^j, \dim(\tilde L)=k\}
\]
where
\[
\tilde\lambda(\tilde L)=\sup\{\sum_{k, l=1}^j c_k \bar c_l Q(u_k,
u_l); \ (c_1, \ldots, c_j)\in \tilde L, \sum_{k=1}^j |c_k|^2=1\}.
\]
Let $L=\{\sum_{l=1}^j c_l u_l; \ (c_1, \ldots, c_j)\in \tilde L\}$.
Then $\tilde\lambda (\tilde L)\ge C_2\lambda(L)$.  Hence
$\tilde\lambda_k\ge C_2\lambda_{k}(Q)$ for all $1\le k\le j$.
Therefore,
\begin{align*}
j\max_{1\le k\le j} Q(u_k, u_k)&\ge \trace(M)=\sum_{k=1}^j Q(u_k,
u_k)=\sum_{k=1}^j\tilde\lambda_k  \ge C_2\sum_{k=1}^j
\lambda_{k}(Q)\\
&\ge C_1 C_2\sum_{k=1}^j k^\eps\ge C_1C_2\int_0^j x^\eps\, dx=C_1
C_2 j^{\eps+1}/(\eps+1).
\end{align*}
Dividing both sides by $j$, we obtain the lemma.
\end{proof}

\begin{proposition}\label{qeigenvalue}
Let $\Omega$ be a smooth bounded pseudoconvex domain in $\C^n$.
Then $\lambda_{j}(\square_q)\le \lambda_{nj}(\square_{q+1})$ for
all $1\le q\le n-1$ and $j$. In particular, if $\square_q$ has
compact resolvent, so is $\square_{q+1}$.
\end{proposition}

\begin{proof}  Let $u=\sumprime_{|J|=q+1} u_Jd\bar z_J\in
C^\infty(\ov{\Omega})\cap \dom(Q_{q+1})$.  Write
\[
u=\frac{1}{(q+1)!}\sum_{|J|=q+1} u_Jd\bar
z_J=\frac{1}{q+1}\sum_{j=1}^n\left(\frac{(-1)^q}{q!}\sum_{|K|=q}
u_{jK}d\bar z_K\right)\wedge d\bar z_j=\frac{1}{q+1}\sum_{j=1}^n
u_j\wedge d\bar z_j
\]
where the $u_j$'s are $(0, q)$-forms defined by the expression in
the parenthesis in the above equalities. It is easy to see that
$u_j\in C^\infty(\ov{\Omega})\cap \dom(Q_q)$ and
$\sum_{j=1}^n\|u_j\|^2=(q+1)\|u\|^2$. Moreover, by the Kohn-Morrey
formula, we have
\begin{align}\label{compare}
\sum_{j=1}^n Q_q(u_j, u_j)&=(q+1)\sumprime_{|J|=q+1}\sum_{j}^n
\int_\Omega \left|\frac{\partial u_J}{\partial \bar z_j}\right|^2
dV+q\sumprime_{|K|=q}\sum_{j, k=1}^n
\int_{b\Omega}\frac{\partial^2\rho(z)}{\partial
z_j\partial \bar z_k} u_{jK}\bar u_{kK} dS\notag \\
&\le (q+1) Q_{q+1}(u, u)
\end{align}
where $\rho$ is any defining function of $\Omega$ whose gradient
has unit length on $b\Omega$. Consider $\widetilde{Q}(u, u)
=\sum_{j=1}^n Q_q(u_j, u_j)$ as a quadratic form on
$\oplus_{j=1}^n L^2_{(0, q)}(\Omega)$. The associated self-adjoint
operator is then $\widetilde\square=\oplus_{j=1}^n \square_q$. Let
$\tilde\lambda_j$ be the number defined by \eqref{minmax} with $Q$
replaced by $\widetilde{Q}$.  We then have
$\lambda_{j}(\square_{q+1})\ge \tilde\lambda_j$. If $\square_q$
has compact resolvent, so does $\widetilde\square$. In this case,
$\tilde\lambda_{nj}=\lambda_{j}(\square_q)$. If $\square_q$ does
not have compact resolvent, let $a$ be the bottom of its essential
spectrum.  If $\lambda_{j}(\square_q)<a$ for all positive integer
$j$, then the $\lambda_{j}(\square_q)$'s are again eigenvalues of
finite multiplicity and
$\tilde\lambda_{nj}=\lambda_{j}(\square_q)$. Otherwise, let $j_0$
be the smallest integer such that $\lambda_{j_0}(\square_q)=a$. In
this case, $\lambda_{j}(\square_q)$, $1\le j < j_0$, are
eigenvalues and $\lambda_{j}(\square_q)=a$ for $j\ge j_0$. Hence
$\tilde\lambda_{nj}=\lambda_{j}(\square_q)$ for $1\le j<j_0$ and
$\tilde\lambda_j=a$ for $j>nj_0$. Therefore, for all cases, we
have $\lambda_{j}(\square_q)=\tilde\lambda_{nj}$.  We thus
conclude the proof of the lemma. \end{proof}

Next we recall some elements of a wavelet construction of
Lemari\'{e} and Meyer (\cite{LemarieMeyer86}; see also
\cite{Daubechies88, HernandezWeiss96}). Let $a(t)\in
C^\infty_0(-\infty, \infty)$ be the cut-off function defined by
$a(t)=\exp(1/(t-1)-2/(2t-1))$ on $(1/2,\ 1)$ and $a(t)=0$
elsewhere.  Let
\[
b(t)= \begin{cases} \left(\int_{1/2}^{1+t} a(s)\,
ds\big/\int^1_{1/2} a(s)\, ds\right)^{1/2}, \quad &\text{if
$t<0$;}\\
& \\
\left(\int_{t}^1 a(s)\, ds\big/\int^1_{1/2} a(s)\,
ds\right)^{1/2}, \quad  &\text{if $t\ge 0$}.
\end{cases}
\]
Then $b(t)$ is a smooth function supported in $[-1/2, \ 1]$ and
satisfying $b(t)\equiv 1$ on $[0, \ 1/2]$ and $b^2(t)+b^2(t-1)\equiv
1$ on $[1/2, \ 1]$. Let $c_0$ be the $L^2$-norm of $b(t)$. Let
\[
g_{k, c}(t)=c_0^{-1} c^{1/2} b(ct) e^{2\pi k c t i}
\]
where $c>0$ and $k\in \Z$.

\begin{lemma}\label{wavelett}
For any given $c>0$, $\{g_{k, c}(t)\}_{k\in\Z}$ is an orthonormal
sequence in $L^2$.
\end{lemma}

\begin{proof}
We provide the proof for completeness.  It is easy to see that the
$L^2$-norm of $g_{k, c}(t)$ is 1. For distinct $k, k'\in \Z$, we
have
\begin{align*}
(g_{k, c}\ g_{k', c})&=c_0^{-2} \int_{-\infty}^\infty b^2(t)
e^{2\pi (k-k')ti}\, dt \\
&=c_0^{-2}\left(\int_{-1/2}^0+\int_0^{1/2}+\int_{1/2}^1\right)
b^2(t) e^{2\pi(k-k')ti}\, dt\\
&=c_0^{-2}\left(\int_{1/2}^1 (b^2(t-1)+b^2(t)) e^{2\pi(k-k')ti}\,
dt+\int_0^{1/2} b^2(t) e^{2\pi(k-k')ti}\, dt\right)\\
&=c_0^{-2}\int_0^1 e^{2\pi (k-k')t i}\,
dt=0.\end{align*}\end{proof}

\section{Special holomorphic coordinates and
non-isotropic bidiscs}\label{polydisc}

The non-isotropic geometry of finite type boundaries in $\C^2$ has
been studied in depth in \cite{NagelSteinWainger85, NRSW89,
Catlin89}. For completeness we provide details below.  The key
difference here is, as noted above, that by using a result of
Forn{\ae}ss and Sibony \cite{FornaessSibony89}, we establish
desirable properties for non-isotropic bidiscs of larger size.

Let $\Omega=\{z\in \C^2 \mid \ r(z)<0\}$ be a smooth bounded
domain with a defining function $r\in C^\infty(\C^2)$. Assume that
$|dr|=1$ on $b\Omega$.  Let
\[
L=\frac{\partial r}{\partial z_2}\frac{\partial}{\partial
z_1}-\frac{\partial r}{\partial z_1}\frac{\partial}{\partial z_2}.
\]
Let $z'\in b\Omega$. For $j, k\ge 1$, let
\[
\scriptl_{jk}\partial\dbar r(z')=\underbrace{L\ldots L}_{
\text{$j-1$ times}}\; \underbrace{\ov{L}\ldots\ov{L}}_{\text{$k-1$
times}}\partial\dbar r(L, \ov{L})(z').
\]
Let $m$ be any positive integer. For any $2\le l\le 2m$, let
\begin{equation}\label{al-def}
A_l(z')=\Big(\sum_{\substack{j+k\le l \\ j, k>0}}
|\scriptl_{jk}\partial\dbar r(z')|^2\Big)^{1/2}.
\end{equation}
Let $\tilde r$ be any defining function for $\Omega$ and let
$\widetilde{L}$ be any non-vanishing complex tangential vector field
of $b\Omega$. Let $\tilde A_l(z')$ be defined by \eqref{al-def} with
$r$ replaced by $\tilde r$ and $L$ by $\widetilde{L}$, then it is
easy to check that  $A_l(z')\approx \tilde A_l(z')$. Thus, whether
or not $A_l(z')$ vanishes is a property that is independent of the
choice of either the defining function or the complex tangential
vector field of $b\Omega$. Furthermore, $b\Omega$ is of finite type
$2m'$ at $z'$ if and only if $A_l(z')$ is 0 for all $2\le l\le
2m'-1$ but is positive for $l=2m'$. For any $\tau>0$, let
\begin{equation}\label{delta-def}
\delta(z', \tau)=\sum_{l=2}^{2m} A_l(z')\tau^l.
\end{equation}
Evidently,
\begin{equation}\label{delta-comp}
\delta(z', \tau)\lesssim \tau^2\quad\text{and}\quad c^{2m}\delta(z',
\tau)\le \delta(z', c\tau)\le c^2\delta(z', \tau),
\end{equation}
for any $\tau$ and $c$ such that $0<\tau, c<1$.  Furthermore,
$b\Omega$ is of finite type $2m$ if and only if $\delta(z',
\tau)\gtrsim \tau^{2m}$ uniformly for all $z'\in b\Omega$ and
$\delta(z'_0, \tau)\lesssim \tau^{2m}$ for some $z'_0\in b\Omega$.

Let $U$ be a neighborhood of a boundary point.  Assume without loss
of generality that $|\partial r/\partial z_2|\gtrsim 1$ on $U$.
After a change of (global) holomorphic coordinates of the form
\begin{equation}\label{coordinate-1}
(\xi_1, \xi_2)=\Phi'(z_1, z_2)=(z_1-z_1',\ 2\frac{\partial
r}{\partial z_2}(z')(z_2-z_2')+\sum_{k=1}^{2m}
\alpha_k(z')(z_1-z_1')^k),
\end{equation}
we have that the image $\Omega'=\Phi'(\Omega)$ is defined by
\begin{equation}\label{defining-0}
\rho(\xi)=r((\Phi')^{-1}(\xi))=\Re\xi_2+\sum_{\substack{2\le
j+k\le 2m
\\j, k>0}} a_{jk}(z')\xi_1^j\bar\xi_1^k+O(|\xi_1|^{2m+1}+
|\xi_2||\xi|).
\end{equation}
The $\alpha_k(z')$'s and $a_{jk}(z')$'s depend smoothly on $z'$, and
they are unique in the sense that if after a holomorphic change of
coordinates $\check\Phi'$ of the form \eqref{coordinate-1} but with
possible different $\alpha_k(z')$'s, $r((\check\Phi')^{-1}(\xi))$ is
in the form of \eqref{defining-0} but with possible different
$a_{jk}(z')$'s, then $\Phi'=\check\Phi'$. (See \cite{Catlin89} for a
detailed discussion on the above coordinates.)  Solving $\Re \xi_2$
in terms of the other variables, we have that $b\Omega'$ is defined
near the origin by $\tilde\rho(\xi)=\Re\xi_2+\tilde h(\xi_1, \Im
\xi_2)$, where
\begin{equation}\label{defining-1}
\tilde h(\xi_1, \Im\xi_2)=\sum_{\substack{2\le j+k\le 2m \\j,
k>0}} \tilde a_{jk}(z')\xi_1^j\bar\xi_1^k+O(|\xi_1|^{2m+1}+
|\Im\xi_2||\xi_1|+|\Im\xi_2|^2).
\end{equation}
It is easy to see that
\begin{equation}\label{al-equiv}
A_l(z')\approx \sum_{\substack{j+k\le l \\ j, k>0}}
|a_{jk}(z')|\approx\sum_{\substack{j+k\le l \\ j, k>0}}|\tilde
a_{jk}(z')|,
\end{equation}
for $2\le l\le 2m$.

Write
\begin{equation}\label{defining-2}
\tilde h(\xi_1,
\Im\xi_2)=\widetilde{P}(\xi_1)+(\Im\xi_2)\widetilde{Q}(\xi_1)+
O(|\xi_1|^{2m+1}+|\Im\xi_2|^2+|\Im\xi_2||\xi_1|^{m+1}),
\end{equation} where
\[
\widetilde{P}(\xi_1)=\sum_{\substack{2\le j+k\le 2m \\j, k>0}}
\tilde a_{jk}(z')\xi_1^j\bar\xi_1^k, \quad
\widetilde{Q}(\xi_1)=\sum_{\substack{1\le j+k\le m}}\tilde
b_{jk}(z')\xi_1^j\bar\xi_1^k.
\]
The harmonic terms can be expunged from the polynomial
$\widetilde{Q}$ without introducing harmonic terms into the
$\widetilde{P}$ term by a change of (local) holomorphic coordinates
of the form
\[
(\tilde\xi_1, \tilde\xi_2)=\Phi^*(\xi_1, \xi_2)=(\xi_1, \xi_2
\prod_{j=1}^{m}(1-\beta_j(z')\xi_1^j)).
\]
(See\cite{FornaessSibony89}.)  Finally, after another change of
coordinates of the form
\[ (\zeta_1, \zeta_2)=\widehat{\Phi}(\tilde\xi_1,
\tilde\xi_2)=(\tilde\xi_1, \tilde\xi_2-\gamma(z')(\tilde{\xi}_2)^2),
\]
we can also eliminate the term $\gamma(z')(\Im\tilde\xi_2)^2$ from
the remainder of the Taylor expansion without introducing harmonic
terms into the $\widetilde{P}$ and $\widetilde{Q}$ terms, and obtain
that $b\Omega$ is defined near $z'$ in the new $(\zeta_1,
\zeta_2)$-coordinates by $\hat\rho=\Re\zeta_2+\hat h(\zeta_1,
\Im\zeta_2)=0$ with
\begin{equation}\label{defining-5}
\hat h(\zeta_1, \Im\zeta_2)=\widehat P(\zeta_1)+
(\Im\zeta_2)\widehat Q(\zeta_1) +O\big(|\zeta_1|^{2m+1}
+|\Im\zeta_2||\zeta_1|^{m+1}+|\Im\zeta_2|^2|\zeta_1|\big),
\end{equation}
where
\[
\widehat
P(\zeta_1)=\sum_{l=2}^{2m}\widehat{P}_l(\zeta_1)=\sum_{l=2}^{2m}
\sum_{\substack{j+k=l
\\j, k>0}} \hat{a}_{jk}(z')\zeta_1^j\bar\zeta_1^k
\]
and
\[
\widehat
Q(\zeta_1)=\sum_{l=2}^{m}\widehat{Q}_l(\zeta_1)=\sum_{l=2}^{m}
\sum_{\substack{j+k=l
\\j, k>0}} \hat{b}_{jk}(z')\zeta_1^j\zeta_1^k.
\]
Let
\[
\widehat{A}_l(z')=\Big(\sum_{\substack{j+k\le l\\ j, k>0}} |\hat
a_{jk}(z')|^2\Big)^{1/2} \quad \text{and}\quad
\widehat{B}_l(z')=\Big(\sum_{\substack{j+k\le l\\ j, k>0}} |\hat
b_{jk}(z')|^2\Big)^{1/2}.
\]
We now summarize what we have obtained from these changes of
holomorphic coordinates. For any $z'\in U\cap b\Omega$, there exists
a neighborhood $U_{z'}$ of $z'$ and a biholomorphic map
$\zeta=\widehat\Psi_{z'}(z)=\widehat\Phi\circ\Phi^*\circ\Phi'(z)$
from $U_{z'}$ onto a ball $B(0, \eps_0)$ of uniform radius $\eps_0$
such that
\begin{enumerate}
\item[(A-1)] $\widehat\Psi_{z'}$ depends smoothly on $z'$ and its
components are holomorphic polynomials of degrees $\le m^2+5m$ for
each $z'$. Moreover, the Jacobian determinant $J\Phi_{z'}$ of
$\Phi_{z'}$ is uniformly bounded from above and below on $U_{z'}$.

\item[(A-2)] $\widehat\Psi_{z'}(z')=0$ and
$\Phi_{z'}(U_{z'}\cap\Omega)=\{\zeta\in B(0, \eps_0) \mid
\hat\rho(\zeta)=\Re \zeta_2+\hat h(\zeta_1, \Im \zeta_2)<0\}$,
where $\hat h(\zeta_1, \Im\zeta_2)$ is in the form of
\eqref{defining-5}.

\item[(A-3)] There exist positive constants $C_1$ and $C_2$
independent of $z'$ such that $C_1A_l(z')\le \widehat{A}_l(z')\le
C_2 A_l(z')$ for $2\le l\le 2m$.
\end{enumerate}

Notice that these properties hold for any smooth bounded domain
$\Omega$.  From now on, we will assume that $\Omega$ is pseudoconvex
of finite type $2m$.  When these assumptions come into play, then it
follows from \cite{FornaessSibony89} that
\begin{equation}\label{FS-1}
\sum_{l=2}^m\|\widehat{Q}_l|\|_\infty|\zeta_1|^l\lesssim
|\zeta_1|\Big(\sum_{l=2}^{2m}\|\widehat{P}_l\|_\infty|\zeta_1|^l
\Big)^{1/2},
\end{equation}
where $\|P\|_\infty$ denotes the sup-norm of $P(\zeta_1)$ on
$|\zeta_1|=1$. It is easy to see that
\[
\|\sum_{j+k=l} c_{jk}\zeta_1^j\bar\zeta_1^k\|_\infty\approx
\sum_{j+k=l}|c_{jk}|.
\]
Therefore, in light of (A-3) and \eqref{FS-1}, we have
\begin{equation}\label{FS-2}
\sum_{l=2}^m \widehat{B}_l(z')\tau^l\lesssim \tau (\delta(z',
\tau))^{1/2}
\end{equation}
for $0<\tau<1$.

Let
\[
\widehat{P}_\tau(z')=\{|\zeta_1|<\tau, |\zeta_2|<\delta(z',
\tau)^{1/2}\},
\]
and let
\[
R_\tau(z')=(\widehat{\Psi}_{z'})^{-1}(\widehat P_\tau(z')).
\]

We now study the non-isotropic ``bidiscs" $R_\tau(z')$. Notice that
the size of $\widehat{P}_\tau(z')$ is different from those used to
study the Bergman kernel in \cite{Catlin89, Mcneal89, NRSW89}. Here
we have $|\zeta_2|<(\delta(z', \tau))^{1/2}$ instead of
$|\zeta_2|<\delta(z', \tau)$.  This seems to be crucial in our
analysis. Let $\tau_0$ be a sufficiently small positive constant
such that $\widehat{P}_{\tau_0}(z')\subset B(0, \eps_0)$ for all
$z'\in b\Omega$. Let
\begin{equation}
\begin{split}\label{class-s}
\widehat{\scripts}^{a, b}_{z'}=&\big\{f\in C^\infty(B(0, \eps_0))
\mid \forall
j, k\ge 0,\ \exists C_{jk} >0\text{ such that }\\
&\quad |D^j_{\zeta_1}D^k_{\zeta_2} f(\zeta)|\le
C_{jk}\tau^{a-j}(\delta(z', \tau))^{b-k/2},  \forall\zeta\in
\widehat{P}_\tau(z'), \  \ \forall\tau\in (0, \tau_0) \big\}.
\end{split}
\end{equation}
(Here $D^j_{\zeta_l}$ denotes the partial derivatives of order $j$
with respect to $\zeta_l$ or $\bar\zeta_l$.)  The following facts
can be checked easily:
\begin{itemize}
\item[(C-1)] If $f\in \widehat{\scripts}^{a, b}_{z'}$, then $D_{\zeta_1}^j
D^k_{\zeta_2} f\in \widehat{\scripts}^{a-j, b-k/2}_{z'}$.

\item[(C-2)] If $f\in \widehat{\scripts}^{a, b}_{z'}$ and $g\in
\widehat{\scripts}^{c, d}_{z'}$, then $fg\in
\widehat{\scripts}^{a+c, b+d}_{z'}$.
\end{itemize}

\begin{lemma}\label{h-estimate} Under the pseudoconvexity and
finite type assumptions on $\Omega$, $\hat h(\zeta_1, \Im\zeta_2)\in
\widehat{\scripts}^{0, 1}_{z'}$.
\end{lemma}
\begin{proof}  Write $\delta=\delta(z', \tau)$ and
let $\hat h_\tau(w_1, w_2)=(1/\delta)\hat h(\tau w_1,
\delta^{1/2}\Im w_2)$. Since $\tau^{2m}\lesssim \delta\lesssim
\tau^2$, it follows that the Taylor expansion of $\hat h_\tau$ at
the origin has the form
\begin{equation}\label{h-hat-tau}
\begin{split}
\hat h_\tau(w_1, w_2)=\sum_{l=2}^{2m}\sum_{j+k=l}&\frac{\tau^l\hat
a_{jk}(z')}{\delta}w_1^j\bar w_1^k+(\Im w_2)\sum_{l=2}^m\sum_{j+k=l}
\frac{\tau^l\hat
b_{jk}(z')}{\delta^{1/2}}w_1^j\bar w_1^k\\
&+O\big(\tau (|w_1|^{2m+1} +|\Im w_2||w_1|^{m+1}+|\Im
w_2|^2)\big).
\end{split}
\end{equation}
Notice that the Taylor coefficients in the first sum above have
modulus $\lesssim 1$ by property (A-3) and those in the second sum
have modulus $\lesssim \tau$ by \eqref{FS-2}. The coefficients of
the Taylor expansion of the remainder are also $\lesssim \tau$ as
shown above.  Therefore, $\hat h(\zeta_1, \Im\zeta_2)\in
\widehat{\scripts}^{0, 1}_{z'}$.
\end{proof}

The next two lemmas establish the doubling and engulfing properties
for the non-isotropic bidiscs $R_\tau(z')$ (cf. \cite{Catlin89}).

\begin{lemma}\label{delta-equiv} Under the same assumptions, if
$z''\in R_\tau(z')\cap b\Omega$, then $\delta(z'', \tau)\approx
\delta(z', \tau)$.
\end{lemma}
\begin{proof}  We shall use the above lemma. (Compare the proof
of Proposition~1.3 in \cite{Catlin89}.) Details are provided for the
reader's convenience. Let $\hat
r(z)=\hat\rho(\widehat{\Psi}^{-1}_{z'}(z))$. Let
\[
L^*=\dfrac{\partial}{\partial\zeta_1}-
\dfrac{\partial\hat\rho}{\partial\zeta_1}\big(\dfrac{\partial\hat\rho}
{\partial\zeta_2}\big)^{-1}\dfrac{\partial}{\partial\zeta_2}
=\dfrac{\partial}{\partial\zeta_1}-\dfrac{\partial\hat h}{\partial
\zeta_1}\big(\dfrac{1}{2}+\dfrac{\partial \hat h}{\partial
\zeta_2}\big)^{-1}\dfrac{\partial}{\partial\zeta_2}.
\]
Notice that the coefficient of $\partial/\partial\zeta_2$ in the
above expression of $L^*$ belongs to $\scripts^{-1, 1}_{z'}$.  Let
$L'=\big(\Psi_{z'}^{-1}\big)_*(L^*)$. Write
\[
\scriptl_{j, k}'\partial\dbar\hat r(z)=\underbrace{L'\ldots
L'}_{j-1\text{ times
}}\;\underbrace{\ov{L'}\ldots\ov{L'}}_{k-1\text{ times
}}\,\partial\bar\partial\hat r(L', \ov{L'})
\]
and
\[
\scriptl_{j,
k}^*\partial\dbar\hat\rho(\zeta)=\underbrace{L^*\ldots
L^*}_{j-1\text{ times
}}\;\underbrace{\ov{L^*}\ldots\ov{L^*}}_{k-1\text{ times
}}\;\partial\bar\partial\hat\rho(L^*, \ov{L^*}).
\]
Then by functoriality, for $z\in U_{z'}$,
\[
\scriptl_{j, k}'\partial\dbar\hat r(z)=\scriptl_{j,
k}^*\partial\dbar\hat\rho(\zeta),
\]
where $\zeta=\widehat\Psi_{z'}(z)$.  It is easy to see that
$\scriptl_{j, k}^*\partial\dbar\hat\rho(\zeta)\in \scripts^{-j-k,
1}_{z'}$; in fact,
\begin{equation}\label{levi-h}
\scriptl_{j, k}^*\partial\dbar\hat\rho(\zeta)=
\frac{\partial^{j+k}\hat
h(\zeta)}{\partial\zeta_1^j\partial\bar\zeta_1^k}+s_{-j-k+1}(\zeta)
\end{equation}
for some $s_{-j-k+1}(\zeta)\in \scripts^{-j-k+1, 1}_{z'}$. It
follows that when $z''\in R_\tau(z')\cap b\Omega$,
\[
\begin{aligned}
A_l(z'')&\approx \max\big\{\big| \scriptl_{j, k}'\partial\dbar\hat
r(z'')\big|; \; 2\le j+k \le l \big\}\\
&=\max\big\{\big| \scriptl_{j, k}^*\partial\dbar\hat\rho
(\zeta'')\big|;\; 2\le j+k \le l \big\}\\
&\lesssim \tau^{-l}\delta(z', \tau).
\end{aligned}
\]
Therefore, $\delta(z'', \tau)\lesssim \delta(z', \tau)$.  We now
prove the estimate in the opposite direction.   From the definition
of $\delta(z', \tau)$ we know that there exist $j_0, k_0>0$ with
$j_0+k_0=l_0\le 2m$ such that
\[
|\scriptl^*_{j_0, k_0}\partial\dbar\hat\rho(0)|\gtrsim
\tau^{-l_0}\delta(z', \tau),
\]
where the constant in the above estimate depending only on $m$. Now
let $z''\in R_{\eps\tau}(z')$ where $\eps$ is a sufficiently small
constant to be determined. By \eqref{levi-h} and \eqref{delta-comp},
we have
\[
\begin{aligned}
\big|\scriptl^*_{j_0, k_0}\partial\dbar\hat\rho(\zeta'')-
\scriptl^*_{j_0, k_0}\partial\dbar\hat\rho(0)\big|&\lesssim
\big|\frac{\partial^{l_0}\hat
h(\zeta'')}{\partial\zeta_1^{j_0}\partial\bar\zeta_1^{k_0}}-
\frac{\partial^{l_0}\hat
h(0)}{\partial\zeta_1^{j_0}\partial\bar\zeta_1^{k_0}}\big|
+\tau^{-l_0+1}\delta(z', \tau)\\
&\lesssim \big(\tau^{-l_0-1}\delta(z',
\tau)\big)|\zeta''_1|+\big(\tau^{-l_0}\delta(z',
\tau)^{1/2}\big)|\zeta''_2|
+\tau^{-l_0+1}\delta(z', \tau)\\
&\lesssim (\eps+\tau) \tau^{-l_0}\delta(z', \tau).\\
\end{aligned}
\]
Therefore when both $\tau$ and $\eps$ are sufficiently small, we
have
\[
\delta(z'', \tau)\gtrsim \tau^{l_0}\big|\scriptl_{j_0,
k_0}\partial\dbar\hat\rho (\zeta'')\big|\gtrsim \delta(z', \tau).
\]
We then conclude the proof by replacing $\eps\tau$ by $\tau$ and
using \eqref{delta-comp}.
\end{proof}

\begin{lemma}\label{engulfing} If $z''\in R_\tau(z')\cap b\Omega$, then there
exists a positive constant $C$ such that
\begin{equation}\label{engulfing-1}
R_\tau(z')\subset R_{C\tau}(z'')\quad \text{and}\quad
R_{\tau}(z'')\subset R_{C\tau}(z').
\end{equation}
\end{lemma}

\begin{proof} It follows from
$\Phi'(R_\tau(z'))=(\widehat{\Phi}\circ\Phi^*)^{-1}(\widehat{P}_\tau(z'))$
that
\[ \{|\xi_1|<C^{-1}\tau, \ |\xi_2|<\delta(z', C^{-1}
\tau)^{1/2}\} \subset \Phi'(R_\tau(z')) \subset \{|\xi_1|<C\tau, \
|\xi_2|<\delta(z', C\tau)^{1/2}\}
\]
for some constant $C>0$.  Thus $\xi''=\Phi'(z'')\in
\{|\xi_1|<C\tau, \ |\xi_2|<\delta(z', C\tau)^{1/2}\}$. After a
change of coordinates of form
\begin{equation}\label{tilde-xi}
(\tilde\xi_1, \tilde\xi_2)=\Psi''(\xi_1, \xi_2)=(\xi_1-\xi_1'',\
2\frac{\partial\rho}{\partial
\xi_2}(\xi'')(\xi_2-\xi_2'')+\sum_{k=1}^{2m}
e_k(z'')(\xi_1-\xi_1'')^k),
\end{equation}
we have
\[
r((\Psi''\circ\Phi')^{-1}(\tilde\xi))=\sum_{\substack{2\le j+k\le
2m
\\j, k>0}}
a_{jk}(z''){\tilde\xi_1}^j{\bar{\tilde\xi}}_1^k+O(|\tilde\xi_1|^{2m+1}+
|\tilde\xi_1||\tilde\xi|).
\]
The $e_k(z'')$'s are determined inductively as follows.
\[
e_1(z'')=2\frac{\partial\rho}{\partial\xi_1}(\xi''), \;
e_k(z'')=\frac{2}{k!}\frac{\partial^k\rho_k}{\partial\xi^k_1}(0), \;
k\ge 2,
\]
where
\[
\rho_1=\rho, \; \rho_k=\rho_{k-1}\circ(\phi_{k-1})^{-1}, \; k\ge 2,
\]
and
\[
\phi_1 =\big(\xi_1-\xi_1'',
2\frac{\partial\rho}{\partial\xi_2}(\xi'')\big(\xi_2-\xi_2'')+
e_1(z'')(\xi_1-\xi_1'')\big)\big), \; \phi_k =\big(\xi_1,
\xi_2+e_k(z'')\xi_1^k\big), \; k\ge 2.
\]
It follows from \eqref{defining-0} and \eqref{al-equiv} that
\[
|D_{\xi_1}^l\rho_1(\xi)|\lesssim \tau^{-l}(\delta(z', \tau))^{1/2},
\qquad\text{for}\quad |\xi_1|\lesssim \tau, \; |\xi_2|\lesssim
(\delta(z', \tau))^{1/2}.
\]
By induction on $k$, we obtain that
\begin{equation}\label{e-k}
|D^l_{\xi_1}\rho_k(\xi)|\lesssim \tau^{-l}(\delta(z', \tau))^{1/2}
\quad\text{and}\quad |e_k(z'')|\lesssim \tau^{-l}(\delta(z',
\tau))^{1/2}.
\end{equation}
By the uniqueness in the sense noted after \eqref{defining-0}, we
obtain as above that
\[
\{|\tilde\xi_1|<C^{-1}\tau, \ |\tilde\xi_2|<\delta(z'',
C^{-1}\tau)^{1/2}\} \subset (\Psi''\circ\Phi')(R_\tau(z''))
\subset \{|\tilde\xi_1|<C\tau, \ |\tilde\xi_2|<\delta(z'',
C\tau)^{1/2}\}.
\]
It follows from Lemma~\ref{delta-equiv}, \eqref{e-k}, and
\eqref{tilde-xi} that if $\xi\in \Phi'(R_\tau(z'))$, then
\[
|\tilde\xi_1|\lesssim \tau\quad\text{and}\quad
|\tilde\xi_2|\lesssim (\delta(z'', \tau))^{1/2}.
\]
Thus, $R_\tau(z')\subset R_{C\tau}(z'')$.  Similarly,
$R_\tau(z'')\subset R_{C\tau}(z')$. \end{proof}

Denote by $d(z)$ the Euclidean distance from $z$ to $b\Omega$. Let
$\pi(z)$ be the projection from a neighborhood of $b\Omega$ onto
$b\Omega$ such that $|z-\pi(z)|=d(z)\approx r(z)$.  Denote by
$\chi_A$ the characteristic function for a set $A$. Let
\[
A_\tau=\{z\in\Omega \mid d(z)<(\delta(\pi(z), \tau))^{1/2}\}.
\]
The following lemma is  an easy consequence of
Lemma~\ref{delta-equiv} and Lemma~\ref{engulfing}.

\begin{lemma}\label{delta-int} For any $\alpha\in \R$, there exists a
sufficiently large constant $C>0$ such that for  any sufficiently
small $\tau>0$ and for any $z\in A_{\tau}$,
\[
\chi_{A_{C^{-1}\tau}}(z)\lesssim \tau^{-2}(\delta(\pi(z),
\tau))^{-\alpha-\frac{1}{2}}
\int_{b\Omega}\chi_{R_\tau(z')\cap\Omega}(z)(\delta(z',
\tau))^{\alpha}\, dS(z')\lesssim \chi_{{A_{C\tau}}}(z).
\]
\end{lemma}

\begin{proof}  It is easy to see that if $z\in R_\tau(z')$ then
$\pi(z)\in R_{C\tau}(z')$.  Thus by Lemma~\ref{delta-equiv},
$\delta(z', \tau)\approx \delta(\pi(z), \tau)$, and by
Lemma~\ref{engulfing}, $z'\in R_{C\tau}(\pi(z))$. It follows that
\[
A_{C^{-1}\tau}\subset \cup_{z'\in b\Omega} R_\tau(z')\cap
\Omega\subset A_{C\tau},
\]
and for any $z\in A_\tau$,
\[
b\Omega\cap R_{C^{-1}\tau}(\pi(z))\subset \{z'\in b\Omega \mid
z\in R_{\tau}(z')\}\subset b\Omega \cap R_{C\tau}(\pi(z)).
\]
Thus
\[
\begin{aligned}
\int_{b\Omega}\chi_{R_\tau(z')\cap\Omega}(z)(\delta(z',
\tau))^{\alpha}\, dS(z')&\lesssim (\delta(\pi(z),
\tau))^{\alpha}{\text{Area}}\left(b\Omega\cap
R_{C\tau}(\pi(z))\right)\chi_{{A_{C\tau}}}(z) \\
&\lesssim \tau^2(\delta(\pi(z),
\tau))^{\alpha+1/2}\chi_{{A_{C\tau}}}(z).
\end{aligned}
\]
The other estimate in Lemma~\ref{delta-int} follows similarly.
\end{proof}

\section{Rescale the $\dbar$-Neumann Laplacian}\label{rescaled}

We will keep the notations from the previous section. Let
$\widehat{\Omega}_{z'}=\widehat{\Psi}_{z'}(\Omega\cap U_{z'})$ and
write
\[
\hat h(\zeta_1,
\Im\zeta_2)=f(\zeta_1)+(\Im\zeta_2)g_1(\zeta_1)+(1/2)(\Im\zeta_2)^2
g_2(\zeta_1)+\sigma_3(\zeta_1, \Im\zeta_2),
\]
where $\sigma_3(\zeta_1, \Im\zeta_2)=O(|\Im\zeta_2|^3)$. Then
\[
f(\zeta_1)=\widehat P(\zeta_1)+O(|\zeta_1|^{2m+1}),\ \
g_1(\zeta_1)=\widehat Q(\zeta_1)+O(|\zeta_1|^{m+1}),\ \
g_2(\zeta_1)=O(|\zeta_1|).
\]
It is evident that $f\in \scripts^{0, 1}$, $g_1\in\scripts^{1,
1/2}$, and $g_2\in\scripts^{1, 0}$.

We flatten the boundary before the rescaling. Let
\[(\eta_1,
\eta_2)=\widetilde{\Phi}_{z'}(\zeta_1, \zeta_2)=(\zeta_1,\
\zeta_2+\hat h(\zeta_1, \Im\zeta_2)-F(\zeta_1, \zeta_2)),
\]
where $F(\zeta_1, \zeta_2)=g_2(\zeta_1)(\Re\zeta_2+\hat h(\zeta_1,
\Im\zeta_2))^2/2+i(g_1(\zeta_1)(\Re\zeta_2)+g_2(\zeta_1)(\Re\zeta_2)
(\Im\zeta_2))$. Of course it is not possible to flatten the boundary
with a holomorphic change of variables: The term $F(\zeta_1,
\zeta_2)$ is added to ensure that
$\partial\eta_2/\partial\bar\zeta_2$ vanishes to a desirable higher
order at the origin.  Note that $F\in \widehat{\scripts}^{1,
1}_{z'}$. Let
\[
\widetilde{P}_\tau(z')=\{|\eta_1|<\tau, \ |\eta_2|<(\delta(z',
\tau))^{1/2}\}.
\]
Let $\widetilde{\scripts}^{a, b}_{z'}$ be the class of smooth
functions in $\eta$ on a neighborhood of the origin defined as in
\eqref{class-s} but with $\zeta$ replaced by $\eta$ and $\widehat
P_\tau(z')$ replaced by $\widetilde{P}_\tau(z')$.

\begin{lemma}\label{phi-tilde} There exists a constant
$C>0$ such that
\[
\widetilde{P}_{C^{-1}\tau}(z')\subset
\widetilde{\Phi}_{z'}(\widehat{P}_\tau(z'))\subset
\widetilde{P}_{C\tau}(z').
\]
\end{lemma}

\begin{proof} The inclusion $\widetilde{\Phi}_{z'}(\widehat{P}_\tau(z'))\subset
\widetilde{P}_{C\tau}(z')$ is evident.  Now if $\eta\in
\widetilde{P}_{\tau}(z')$, then
\begin{equation}
\begin{aligned}
|\zeta_2|&=|\eta_2-\hat h(\eta_1, \Im\zeta_2)-F(\eta_1,
\zeta_2)|\\
&\le \delta^{1/2}+|\hat h(\eta_1, \Im\zeta_2)|+|F(\eta_1,
\zeta_2)|\\
&\lesssim \delta^{1/2}+ \tau\delta^{1/2}|\zeta_2|+\tau
|\zeta_2|^2+|\zeta_2|^3.
\end{aligned}
\end{equation}
Thus $|\zeta_2|\lesssim \delta^{1/2}$.  The other inclusion then
follows. \end{proof}

Let $\check\rho(\zeta)=\hat\rho(\zeta)-
(1/2)g_2(\zeta_1)(\hat\rho(\zeta))^2$. Then $\check\rho(\zeta)$ is a
defining function for $b\widehat\Omega_{z'}$ near the origin. Let
$\check r(z)=\check\rho(\widehat{\Psi}_{z'}(z))$. Then $\check r(z)$
is a defining function for $b\Omega\cap U_{z'}$ (shrinking $U_{z'}$
if necessary). Let
\begin{equation}\label{vector-basis}
L_1=\frac{\partial\check r}{\partial z_2}\frac{\partial}{\partial
z_1}-\frac{\partial\check r}{\partial z_1}\frac{\partial}{\partial
z_2}\quad\text{and}\quad L_2=\frac{\partial\check r}{\partial \bar
z_1}\frac{\partial}{\partial z_1}+\frac{\partial\check r}{\partial
\bar z_2}\frac{\partial}{\partial z_2},
\end{equation}
and let
\begin{equation}\label{form-basis}
\omega_1=\frac{\partial\check r}{\partial \bar
z_2}dz_1-\frac{\partial\check r}{\partial\bar z_1}dz_2
\quad\text{and}\quad \omega_2=\frac{\partial\check r}{\partial z_1}
dz_1+\frac{\partial\check r}{\partial z_2}dz_2.
\end{equation}
Then $\{ L_1, L_2\}$ forms an orthogonal basis for $T^{1,0}(\C^2)$
and $\{\omega_1, \omega_2\}$ for $\Lambda^{1, 0}(\C^2)$ on $U_{z'}$.
Denote by $\widehat{L}_1$, $\widehat{L}_2$, $\widehat{\omega}_1$,
and $\widehat{\omega}_2$ the vectors and forms defined as above by
replacing $\check r$ by $\check\rho$, and $z_1$, $z_2$ by $\zeta_1$,
$\zeta_2$ respectively.  Let
$\widetilde{L}_k=(\widetilde{\Phi}_{z'})_*(\widehat{L}_k)$, $k=1,
2$. Write $\zeta_2=\hat s+i\hat t$ and $\eta_2=\tilde s+i\tilde t$.

\begin{lemma}\label{l-tilde} With above notations,
\[
\ov{\widetilde{L}}_1=\big(\frac{1}{2}+\alpha_1\big)\frac{\partial}
{\partial\bar\eta_1}+\big(-\frac{i}{2}\frac{\partial\hat
h}{\partial \bar\zeta_1}+\beta_1\big)\frac{\partial}{\partial
\tilde t}, \ \
\ov{\widetilde{L}}_2=\big(\frac{1}{2}+\alpha_2\big)\frac{\partial}
{\partial\bar\eta_2}+\alpha_3\frac{\partial}{\partial\bar\eta_1}+
\beta_2\frac{\partial}{\partial\tilde
s}+\beta_3\frac{\partial}{\partial \tilde t},
\]
where the $\alpha$'s are in $\widetilde\scripts^{0, 1/2}_{z'}$ and
the $\beta$'s in $\widetilde\scripts^{0, 1}_{z'}$.
\end{lemma}

\begin{proof} By direct computations, we have
\[
\begin{aligned}
\ov{\widetilde{L}}_1&=\frac{\partial\check\rho}{\partial\bar
\zeta_2}\frac{\partial}{\partial\bar\eta_1}+\left(\frac{\partial\check\rho}
{\partial\bar\zeta_2}\frac{\partial\tilde
t}{\partial\bar\zeta_1}-\frac{\partial\check\rho}
{\partial\bar\zeta_1}\frac{\partial\tilde
t}{\partial\bar\zeta_2}\right)\frac{\partial}{\partial\tilde t},\\
\ov{\widetilde{L}}_2&=\left[\frac{\partial\check\rho}{\partial
\zeta_2}\left(1+\frac{\partial\hat h}{\partial\bar
\zeta_2}-\frac{\partial\bar
F}{\partial\bar\zeta_2}\right)+\frac{\partial\check\rho}{\partial
\zeta_1}\left(\frac{\partial\hat h}{\partial\bar
\zeta_1}-\frac{\partial\bar
F}{\partial\bar\zeta_1}\right)\right]\frac{\partial}{\partial\bar\eta_2}+
\frac{\partial\check\rho}{\partial\zeta_1}\frac{\partial}
{\partial\bar\eta_1}\\
&\qquad\qquad+ \left[\frac{\partial\check\rho}{\partial
\zeta_1}\left(\frac{\partial\hat h}{\partial\bar
\zeta_1}-\frac{\partial
F}{\partial\bar\zeta_1}\right)+\frac{\partial\check\rho}{\partial
\zeta_2}\left(\frac{\partial\hat h}{\partial\bar
\zeta_2}-\frac{\partial F}{\partial\bar\zeta_2}\right)\right]
\frac{\partial}{\partial\eta_2}.
\end{aligned}
\]
Note that $\widetilde{\Phi}^*_{z'}$ is an isomorphism from
$\widetilde\scripts^{a, b}_{z'}$ onto $\widehat{\scripts}^{a,
b}_{z'}$. The lemma then follows from the facts that $\hat
h\in\widehat{\scripts}^{0, 1}_{z'}$, $\check\rho-\Re\zeta_2\in
\widehat{\scripts}^{0, 1}_{z'}$,  $F\in \widehat{\scripts}^{1,
1}_{z'}$, and
\[
\frac{\partial\hat h}{\partial\bar\zeta_2}-\frac{\partial F}
{\partial\bar\zeta_2}=-\frac{i}{2}g_2(\hat s+\hat
h)\frac{\partial\hat h}{\partial \hat t}-\frac{1}{2}g_2\hat
h+O(|\hat t|^2)\in \widehat{\scripts}^{0, 1}_{z'}.
\]\end{proof}

We now proceed with the rescaling.  For any positive $\tau$, let
$\delta=\delta(z', \tau)$ be defined by \eqref{delta-def}. Let
\[
(w_1, w_2)=D_{z', \tau}(\eta_1, \eta_2)=(\eta_1/\tau,
\eta_2/\delta).
\]
Let  $\widetilde{\Psi}_{z', \tau}=D_{z',
\tau}\circ\widetilde{\Phi}_{z'}$ and let $\widetilde\Omega_{z',
\tau}=\widetilde{\Psi}_{z', \tau}(\widehat\Omega_{z'})$. (In what
follows, we sometimes suppress the subscript $z'$ for economy of
notations when there is no confusion.) Let
\[
P_\tau(z')=\{|w_1|<1, \ |w_2|<\delta^{-1/2}\}.
\]
Let $\scripts^{a, b}_{z'}$ be the class of functions $f$ depending
smoothly on $w\in \C^2$ and $\tau>0$ such that for any $j, k\ge 0$,
there exists constants $C_{jk}>0$, independent of $\tau$, such that
\[
|D^j_{w_1}D^k_{w_2} f|\le C_{jk} \tau^a\delta^{b+k/2}
\]
on $P_\tau (z')$ for sufficiently small $\tau>0$. Here, as before,
$D^j_{w_l}$ denotes the partial derivatives of order $j$ with
respect to $w_l$ or $\bar w_l$. Clearly, if
$g\in\widetilde\scripts^{a, b}_{z'}$, then $(D^{-1}_\tau)^*(g)\in
\scripts^{a, b}_{z'}$.

Write $w_1=x+iy$ and $w_2=s+it$. Let
\[
\begin{aligned}
\ov{\widetilde{L}}_{1, \tau}&=\tau
D_{\tau*}(\ov{\widetilde{L}}_{1})=\big(\frac{1}{2}+\alpha_1\big)\frac{\partial}
{\partial\bar w_1}+\big(-\frac{i}{2}\frac{\partial\hat h}{\partial
\bar\zeta_1}+\beta_1\big)\frac{\tau}{\delta}\frac{\partial}{\partial t},\\
\ov{\widetilde{L}}_{2, \tau}&=\tau
D_{\tau*}(\ov{\widetilde{L}}_2)=\big(\frac{1}{2}+\alpha_2\big)\frac{\tau}{\delta}\frac{\partial}
{\partial\bar w_2}+\alpha_3\frac{\partial}{\partial\bar w_1}+
\beta_2\frac{\tau}{\delta}\frac{\partial}{\partial
s}+\beta_3\frac{\tau}{\delta}\frac{\partial}{\partial t}.
\end{aligned}
\]
Write $L^0=\widetilde{L}_{1, \tau}$ and
$L^1=\ov{\widetilde{L}}_{1, \tau}$.  For any tuple $(i_1\ldots
i_l)$ of 0's and 1's, define $L^{(i_1\ldots i_l)}$ inductively by
\[
L^{(i_1\ldots i_l)}=[L^{i_l},  L^{(i_1\ldots i_{l-1})}].
\]
Write
\[
L^{(i_1\ldots i_l)}=\lambda^{i_1\ldots
i_l}\frac{\partial}{\partial t}+ a^{i_1\ldots i_l}L^0+b^{i_1\ldots
i_l}L^1.
\]

\begin{lemma}\label{commutator} With the above notations,
\begin{enumerate}
\item $\lambda^{i_1\ldots i_l}\in \scripts^{0, 0}_{z'}$,
$a^{i_1\ldots i_l}\in \scripts^{1, 0}_{z'}$, and $b^{i_1\ldots
i_l}\in\scripts^{1, 0}_{z'}$.

\item $(L^{(i_1\ldots i_l)})^*=-L^{((1-i_1)\ldots
(1-i_l))}+\sigma$ for some $\sigma\in \scripts^{1, 0}_{z'}$.

\item There exists a tuple $(i_1\ldots i_{l_0})$ of length $l_0\le
2m$ such that $|\lambda^{i_1\ldots i_{l_0}}|\gtrsim 1$ on
$P_\tau(z')$.
\end{enumerate}
\end{lemma}

\begin{proof} A direction calculation yields that
\[
\begin{aligned}
L^{(10)}&=\left(-\frac{i}{4}\left(\frac{\partial}{\partial
w_1}\left(\frac{\partial\hat h}{\partial\bar\zeta_1}\right)
+\frac{\partial}{\partial \bar w_1}\left(\frac{\partial\hat
h}{\partial\zeta_1}\right)\right)\frac{\tau}{\delta}+\sigma
\right)\frac{\partial}{\partial t}+a^{10}L^0+b^{10}L^1\\
&=\left(-\frac{i}{2}\frac{\partial^2\hat
h}{\partial\zeta_1\partial\bar\zeta_1}\frac{\tau^2}{\delta}+\sigma\right)
\frac{\partial}{\partial t}+a^{10}L^0+b^{10}L^1
\end{aligned}
\]
with $a^{10}, b^{10}, \sigma\in \scripts^{1, 0}_{z'}$. (Here and in
what follows, $\sigma$ could be different in different appearances,
but is always in $\scripts^{1, 0}_{z'}$.)  It is also easy to see
that
\[
[L^{i_{l+1}}, \frac{\partial}{\partial
t}]=\sigma\frac{\partial}{\partial t} \mod (L^0, L^1)
\]
where the modulus is with coefficients in $\scripts^{1, 0}_{z'}$.
Thus,
\[
L^{(i_1\ldots i_{l+1})}=\left(L^{i_{l+1}}\lambda^{i_1\ldots
i_l}+\lambda^{0i_{l+1}}a^{i_1\ldots
i_l}+\lambda^{1i_{l+1}}b^{i_1\ldots i_l}+\sigma\lambda^{i_1\ldots
i_l}\right)\frac{\partial}{\partial t} \mod (L^0, L^1).
\]
Properties (1) and (2) in the lemma then follow from an easy
inductive argument on $l$.  To prove (3), one notices from the above
formulas that
\[
\lambda^{10i_3\ldots i_l}=-\frac{i}{2^{l+1}}\frac{\partial^l\hat
h}{\partial\zeta^j_1\partial\zeta^k_1}\frac{\tau^l}{\delta}+\sigma,
\]
where $j$ and $k$ are the numbers of the 0's and 1's in
$(10i_3\ldots i_l)$ respectively. It follows from the proof of
Lemma~\ref{delta-equiv} that there exists $j_0, k_0>0$ with
$j_0+k_0=l_0\le 2m$ such that
\[
\left|\frac{\partial^{l_0}\hat h}{\partial \zeta^{j_0}_1\partial
\bar\zeta^{k_0}_1}\frac{\tau^{l_0}}{\delta}\right|\approx 1
\]
on $R_\tau(z')$.  This then implies the last part of the lemma.
\end{proof}

We now define the rescaled $\dbar$-Neumann Laplacian.  Let
$\widehat\scriptg_{\tau}\colon
(L^2(\widetilde{\Omega}_{\tau}))^2\to L^2_{(0,
1)}(\widehat{\Omega})$ be the transformation defined by
\[
\widehat\scriptg_{\tau} (u_1, u_2)=|\det
d\widetilde{\Psi}_\tau|^{1/2}\big(u_1
(\widetilde{\Psi}_\tau)\ov{\widehat\omega}_1+
u_2(\widetilde{\Psi}_\tau)\ov{\widehat \omega}_2\big),
\]
and let $\widehat\scriptg\colon L^2_{(0, 1)}(\widehat\Omega)\to
L^2_{(0, 1)}(\Omega\cap U)$ be defined likewise by
\[
\widehat\scriptg(u_1 \ov{\widehat \omega}_1+u_2 \ov{\widehat
\omega}_2)=J\widehat\Psi(z) \big(u_1
(\widehat\Psi)\ov{\omega}_1+u_2(\widehat\Psi)\ov{\omega}_2\big),
\]
where $J\Psi$ is the Jacobian determinant of $\Psi$. Evidently,
$\|\widehat\scriptg_\tau u\|_{\widehat\Omega}\approx
\|u\|_{\widetilde\Omega_\tau}$ and $\|\widehat\scriptg
u\|_{\widehat\Omega}\approx \|u\|_{\Omega\cap U}$.  Let
$\scriptg_{\tau}=\widehat\scriptg\circ\widehat\scriptg_{\tau}$ and
let $\widetilde\scriptg_{\tau}(u)=(\delta\tau)^{-2} u\circ D_{\tau}$
be the unitary transformation on $L^2$-spaces associated with the
dilation $D_{\tau}$. Let
\begin{equation}\label{Q-tau-def}
Q_{\tau}(u, v)=\tau^2 Q(\scriptg_{\tau} u, \scriptg_{\tau} v)
\end{equation}
be the densely defined, closed sesquilinear form on
$(L^2(\widetilde\Omega_{\tau}))^2$ with
$\dom(Q_{\tau})=\{\scriptg^{-1}_{\tau}(u); \ u\in\dom(Q), \ \supp
u\subset \ov{\Omega}\cap U\}$.  Here $Q(\cdot, \cdot)$ is the
sesquilinear form associated with the $\dbar$-Neumann Laplacian on
$L^2_{(0, 1)}(\Omega)$.

\begin{lemma}\label{Q-equiv} For any $u\in \dom(Q)\cap
C^\infty_0(\ov{\Omega}\cap U)$,
\[
Q(u, u)\approx \widehat{Q} (\widehat{\scriptg}^{-1} u,
\widehat{\scriptg}^{-1} u),
\]
where $\widehat{Q}$ is the sesquilinear form associated with the
$\dbar$-Neumann Laplacian on $L^2_{(0, 1)}(\widehat\Omega)$.
\end{lemma}
\begin{proof} From \cite{Kohn72}, we know that

\begin{equation}\label{Kohn-Morrey-1}
Q(u, u)\approx \|u\|^2_{\Omega}+ \|\ov{L}_1 u\|^2_{\Omega}+
\|\ov{L}_2 u\|^2_{\Omega}+\int_{b\Omega} \big(\partial\dbar\hat
r(L_1, \ov{L}_1)(z)\big)|u|^2\, dS(z).
\end{equation}
It follows from (A-1) that
\[
\begin{aligned}
Q(u, u)&\approx \| \widehat \scriptg^{-1}u\|^2_{\widehat\Omega}+
\|\ov{\widehat L}_1\widehat \scriptg^{-1} u\|^2_{\widehat\Omega}+
\|\widehat\Psi_*(\ov{L}_2)\widehat \scriptg^{-1}
u\|^2_{\widehat\Omega}+\int_{b\widehat\Omega} \big(\partial\dbar\hat
\rho(\widehat L_1, \ov{\widehat
L}_1)(\zeta)\big)|\widehat \scriptg^{-1} u|^2\, dS(\zeta) \\
&\approx \| \widehat \scriptg^{-1}u\|^2_{\widehat\Omega}+
\|\ov{\widehat L}_1\widehat \scriptg^{-1} u\|^2_{\widehat\Omega}+
\|\ov{\widehat L}_2\widehat \scriptg^{-1}
u\|^2_{\widehat\Omega}+\int_{b\widehat\Omega}
\big(\partial\dbar\hat \rho(\widehat L_1, \ov{\widehat
L}_1)(\zeta)\big)|\widehat \scriptg^{-1} u|^2\, dS(\zeta).
\end{aligned}
\]
Thus $Q(u, u)\approx \widehat Q(\widehat \scriptg^{-1} u, \widehat
\scriptg^{-1} u)$. \end{proof}

Let $\widetilde u(\xi', s)=(\scriptf_{\tan}u)(\xi', s)$ be the
tangential Fourier transform of $u$ in the $x'=(x, y, t)$ variables.
Recall that the the tangential Laplacian $\Lambda^s$ is defined by
\[
\scriptf_{\tan}(\Lambda^s u)(\xi',
s)=(1+|\xi'|^2)^{s/2}\widetilde{u}(\xi', s)
\]
and the tangential $L^2$-Sobolev norm of order $s$ by
\[
\vvv u \vvv_s^2=\int_{-\infty}^0 \int_{\R^3} (1+|\xi|^2)^s
|\widetilde{u}(\xi', s)|^2\, d\xi' ds.
\]

\begin{lemma}\label{Q-estimate}  There exists an $\eps>0$ such that
for any sufficiently small $\tau>0$,
\begin{equation}\label{Q-estimate-2}
Q_{\tau}(u, u)\gtrsim \vvv u \vvv^2_\eps+\tau^2\delta^{-2} \vvv
\frac{\partial u}{\partial \bar w_2}\vvv^2_{-1+\eps},
\end{equation}
for all $u\in\dom(Q_\tau)\cap C^\infty_0(P_\tau(z'))$.
\end{lemma}

\begin{proof} By Lemma~\eqref{Q-equiv},
\begin{equation}\label{Q-estimate-3}
\begin{aligned}
Q_\tau(u, u) &\approx \tau^2 \widehat Q (\widehat\scriptg_\tau u,
\widehat\scriptg_\tau u)\\
&\approx
\tau^2\big(\|\widehat{\scriptg}_{\tau}u\|^2_{\widehat\Omega}+
\|\ov{\widehat{L}}_1
\widehat{\scriptg}_{\tau}u\|^2_{\widehat\Omega}+ \|\widehat{L}_1
\widehat{\scriptg}_{\tau}u\|^2_{\widehat\Omega}
+\|\ov{\widehat{L}}_2
\widehat{\scriptg}_{\tau}u\|^2_{\widehat\Omega}\big)\\
&\approx \tau^2\|u\|^2_{\widetilde\Omega_\tau}+
\|\ov{\widetilde{L}}_{1, \tau} u\|^2_{\widetilde\Omega_\tau}+
\|\widetilde{L}_{1, \tau} u\|^2_{\widetilde\Omega_\tau}
+\|\ov{\widetilde{L}}_{2, \tau} u\|^2_{\widetilde\Omega_\tau}.
\end{aligned}
\end{equation}
We first prove that there exists an $\eps>0$ such that
\begin{equation}\label{kohn}
\vvv \frac{\partial u}{\partial t}\vvv^2_{-1+\eps}\lesssim
Q_\tau(u, u).
\end{equation}
This is a direct consequence of Kohn's method \cite{Kohn72}, in
light of \eqref{Q-estimate-3} and Lemma~\ref{commutator}. Since we
need to keep track that the constant in \eqref{kohn} is independent
of $\tau$, we sketch the proof for completeness.  By
Lemma~\ref{commutator} (1) and (3), we have
\begin{equation}\label{kohn1}
\vvv \frac{\partial u}{\partial t}\vvv^2_{\eps-1}\lesssim
\vvv\lambda^{i_1\ldots i_{l_0}}\frac{\partial u}{\partial
t}\vvv^2_{\eps-1}\lesssim \vvv L^{(i_1\ldots i_{l_0})}
u\vvv^2_{\eps-1}+Q_\tau(u, u).
\end{equation}
It remains to estimate $\vvv L^{(i_1\ldots i_{l_0})}
u\vvv^2_{\eps-1}$, which equals
\begin{equation}\label{two-term}
\langle L^{(i_1\ldots i_{l_0-1})} u,
(L^{i_{l_0}})^*\Lambda^{2(\eps-1)}L^{(i_1\ldots i_{l_0})} u\rangle
-\langle L^{i_{l_0}} u, (L^{(i_1\ldots
i_{l_0-1})})^*\Lambda^{2(\eps-1)}L^{(i_1\ldots i_{l_0})} u\rangle.
\end{equation}
The first term above equals
\[
\begin{split}
\langle L^{(i_1\ldots i_{l_0-1})} u,\
&\Lambda^{2(\eps-1)}L^{(i_1\ldots i_{l_0})}(L^{i_{l_0}})^*
u\rangle+ \langle L^{(i_1\ldots i_{l_0-1})} u, \ [(L^{i_{l_0}})^*,
\Lambda^{2(\eps-1)}] L^{(i_1\ldots i_{l_0})} u\rangle \\
&+\langle L^{(i_1\ldots i_{l_0-1})} u, \
\Lambda^{2(\eps-1)}[(L^{i_{l_0}})^*, L^{(i_1\ldots i_{l_0})}]
u\rangle =I+II+III.
\end{split}
\]
We have
\[
\begin{aligned}
|I|&=|\langle (\Lambda^{2(\eps-1)}L^{(i_1\ldots i_{l_0})})^*
L^{(i_1\ldots i_{l_0-1})} u,\ (-L^{(1-i_{l_0})} +\sigma)
u\rangle|\\
&\lesssim \vvv L^{(i_1\ldots i_{l_0-1})}
u\vvv^2_{2\eps-1}+Q_\tau(u, u),
\end{aligned}
\] because
$(\Lambda^{2(\eps-1)}L^{(i_1\ldots i_{l_0})})^*$ is a tangential
pseudodifferntial operator of order $2\eps-1$. Also, since
$[(L^{i_{l_0}})^*, \Lambda^{2(\eps-1)}]$ is of order $2(\eps-1)$,
we have
\[
|II|\lesssim C\vvv L^{(i_1\ldots i_{l_0-1})}
u\vvv^2_{\eps-1}+(1/C)\vvv L^{(i_1\ldots i_{l_0})}
u\vvv^2_{\eps-1}.
\]
Furthermore, as for \eqref{kohn1}, it follows from
Lemma~\ref{commutator} that
\[
\begin{aligned}
|III|&=|\langle L^{(i_1\ldots i_{l_0-1})} u, \
\Lambda^{2(\eps-1)}[-L^{(1-i_{l_0})}+\sigma, L^{(i_1\ldots
i_{l_0})}] u\rangle|\\
&\lesssim C\vvv L^{(i_1\ldots i_{l_0-1})}
u\vvv^2_{\eps-1}+(1/C)\vvv L^{(i_1\ldots i_{l_0})}
u\vvv^2_{\eps-1}+Q_\tau(u, u).
\end{aligned}
\]
The second term in \eqref{two-term} is estimated similarly and is
left to the reader. From these estimates, we then have
\[
\vvv L^{(i_1\ldots i_{l_0})} u\vvv^2_{\eps-1}\lesssim \vvv
L^{(i_1\ldots i_{l_0-1})} u\vvv^2_{2\eps-1}+Q_\tau(u, u).
\]
Let $\eps=2^{-2m}$. Repeating the above arguments, we then obtain
\eqref{kohn}.

Since $\partial/\partial s$ is a linear combination of $L^0$, $L^1$,
and $\ov{\widetilde{L}}_2$ with coefficients in $\scripts^{1,
0}_{z'}$, it follows from \eqref{kohn} that
\begin{equation}\label{kohn-s}
\vvv \frac{\partial u}{\partial s}\vvv^2_{-1+\eps}\lesssim
Q_\tau(u, u).
\end{equation}
Combining \eqref{Q-estimate-3}, \eqref{kohn}, \eqref{kohn-s}, and
Lemma~\ref{l-tilde}, we then have
\[
\vvv\nabla u\vvv^2_{\eps-1}+\tau^2\delta^{-2}\vvv\frac{\partial
u}{\partial\bar w_2}\vvv^2_{\eps-1}\lesssim Q_\tau(u, u).
\]
By applying the Poincar\'{e} inequality to $\widetilde{u}(\xi',
\cdot)$, we know that the left-hand side above dominates
$\|u\|^2$. We thus conclude the proof of the lemma.
\end{proof}

\section{Auxiliary estimates}\label{aux}

For any $\eps$ such that $0<\eps\le 1/2$ and any $\delta>0$,  let
$W_{\eps, \delta}$ be the space of all $u\in L^2(\C^2_{-})$ such
that
\begin{equation}\label{m-tau-norm}
\|u\|^2_{\eps, \delta}= \vvv u \vvv^2_{\eps}+\delta^{-1}\vvv
\frac{\partial u}{\partial\bar w_2}\vvv^2_{-1+\eps} <\infty.
\end{equation}
Let ${\widetilde Q}_{\eps, \delta}$ be the sesquilinear form on
$L^2(\C^2_-)$ associated with the above norm with $\dom({\widetilde
Q}_{\eps, \delta})$ $=W_{\eps,\delta}$.  Let
${\widetilde\square}_{\eps, \delta}$ be the associated densely
defined, self-adjoint operator on $L^2(\C^2_-)$ and let
$\widetilde{N}_{\eps, \delta}$ be its inverse. Let $\chi (w_1, w_2)$
be a smooth cut-off function supported on $\{|w_1|<1, |w_2|<1\}$ and
identically 1 on $\{|w_1|<1/2, |w_2|<1/2\}$.  Let $\chi_\delta(w_1,
w_2)=\chi(w_1, \delta^{1/2} w_2)$.  We now study the spectral
behavior of $\chi_\delta\widetilde{N}_{\eps, \delta}$ as $\delta\to
0^+$. Let $P$ be the orthogonal projection from $L^2(\C^2_-)$ onto
$H=\{u\in L^2(\C^2_-) \mid \partial u/\partial\bar w_2=0\}$. Namely,
$P$ is the partial Bergman projection in the $w_2$-variable.

\begin{lemma}\label{Bergman}  For all $\delta>0$ and $u\in W_{\eps, \delta}$,
\begin{equation}\label{Bergman-estimate}
\|(I-P)u\|_{\eps, \delta}\lesssim \|u\|_{\eps, \delta}; \ \
\vvv\frac{\partial}{\partial s}(I-P)u\vvv_{-1+\eps}+\vvv
\frac{\partial}{\partial t}(I-P)u\vvv_{-1+\eps} \approx \vvv
\frac{\partial u}{\partial \bar w_2}\vvv_{-1+\eps}.
\end{equation}
\end{lemma}

\begin{proof}  This lemma follows from standard elliptic
theory (cf. \cite{Metivier81}).  We provide the proof for
completeness. Recall that $I-P=4\frac{\partial}{\partial
w_2}G\frac{\partial}{\partial\bar w_2}$, where $G$ is the Green's
operator in $w_2$-variable ({\it i.e.}, the inverse of
$-\Delta_{w_2}$). Throughout this section, we will use
$\zeta_1=\xi+i\eta$ and $\zeta_2=\mu+i\nu$ to denote the dual
variables of $w_1=x+iy$ and $w_2=s+it$ in the Fourier transform.
Recall that $\widetilde u$ denotes the tangential Fourier transform
of $u$ in the $(x, y, t)$ variables.  We have
\[
\vvv (I-P)u\vvv^2_{\eps}=\int_{-\infty}^0\,
ds\int_{\R^3}(1+|\zeta_1|^2+\nu^2)^{\eps}
|2(\frac{\partial}{\partial s}+\nu)\widetilde{G(\frac{\partial
u}{\partial \bar w_2})}|^2\,d\xi d\eta dt.
\]
It is easy to check that $\widetilde{G(u)}=-E_+E_-\widetilde{u}$,
where
\[
(E_-\widetilde{u})(\zeta_1, s, \nu)=\int_{-\infty}^s
e^{-|\nu|(s-s')}\widetilde{u}(\zeta_1, s', \nu)\, ds'
\]
and
\[
(E_+\widetilde{u})(\zeta_1, s, \nu)=-\int_s^0 e^{|\nu|(s-s')}
\widetilde{u}(\zeta_1, s', \nu)\, ds'.
\]
(See, {\it e.g.}, Chapter III in \cite{Treves75}.)  Using the
identities $\frac{\partial (E_+\widetilde{u})}{\partial
s}=|\nu|E_+\widetilde{u}+\widetilde{u}$ and $E_-(\frac{\partial
\widetilde{u}}{\partial
s})=-|\nu|E_-(\widetilde{u})+\widetilde{u}$, we obtain that
\[
2(\frac{\partial}{\partial s}+\nu)\widetilde{G(\frac{\partial
u}{\partial \bar w_2})}=(\nu+|\nu|)^2E_+E_-\widetilde{u}-
(\nu+|\nu|)(E_+{\widetilde{u}}-E_-{\widetilde{u}})-\widetilde{u}.
\]
Since by the Minkowski inequality,
\[\int_{-\infty}^0 |E\widetilde{u}(\zeta_1, s, \nu)|^2\, ds\le
|\nu|^{-2}\int_{-\infty}^0 |\widetilde{u}(\zeta_1, s, \nu)|^2\, ds
\]
holds for both $E_+$ and $E_-$, we obtain that $\vvv
(I-P)u\vvv_{\eps}\lesssim \vvv u\vvv_{\eps}$. The first inequality
then follows.  The second inequality is treated similarly and its
proof is left to the reader. \end{proof}

\begin{lemma}\label{eigen-estimate-1} For sufficiently small $\delta>0$
and sufficiently large $j$,
\[
\lambda_j(\chi_\delta\widetilde{N}^{1/2}_{\eps, \delta})\lesssim
(1+j\delta^{1/2})^{-\eps/4}.
\]
\end{lemma}

\begin{proof} For $u\in L^2(\C^+_-)$, we write
$\widetilde{N}^{1/2}_{\eps,
\delta}u=(I-P)\widetilde{N}^{1/2}_{\eps,
\delta}u+P\widetilde{N}^{1/2}_{\eps, \delta}u=v_1+v_2$.  We first
study $\chi_\delta(I-P)\widetilde{N}_{\eps, \delta}$.  We extend
$v_1$ evenly to $s>0$ by letting $v_1(w_1, s+it)=v_1(w_1, -s+it)$.
Denote by $\widehat v_1$ the Fourier transform of the extended
$v_1$ in all variables. Then by Lemma~\ref{Bergman},
\begin{equation}\label{v-1-estimate}
\begin{aligned}
\|\widetilde{N}^{1/2}_{\eps, \delta} u\|^2_{\eps, \delta}&\gtrsim
\int_{\R^4}\left( (1+|\zeta_1|^2+\nu^2)^\eps
+\delta^{-1}|\zeta_2|^2(1+|\zeta_1|^2+\nu^2)^{-1+\eps}\right)
|\widehat{v}_1|^2\, dV(\zeta)\\
&\gtrsim \int_{\R^4} (1+|\zeta_1|^2+\delta^{-1}|\zeta_2|^2)^\eps
|\widehat{v}_1|^2\, dV(\zeta) \equiv \|v_1\|^{'2}_{\eps, \delta},
\end{aligned}
\end{equation}
where in the last estimate we use the following simple inequality:
$a^\eps b^{1-\eps}\le \eps a+(1-\eps) b$. Let $\Delta'_{\eps,
\delta}$ be the Dirichlet realization of the self-adjoint operator
associated with the sesquilinear form that defines the norm
$\|\cdot \|'_{\eps, \delta}$ on $\{|w_1|<1, \
|w_2|<\delta^{-1/2}\}$.  As such, we have $\|v\|'_{\eps,
\delta}=\|(\Delta'_{\eps, \delta})^{1/2}v\|$. Let $S_{\delta} u
(w_1, w_2)=\delta^{1/2} u(w_1, \delta^{1/2}w_2)$. Then $S_\delta$
is an isometry on $L^2(\C^2_-)$. Furthermore, it is easy to see
that $\|v\|^{'}_{\eps, \delta}=\|S_\delta v\|^{'}_{\eps, 1}$.  It
follows that
\begin{equation}\label{weyl}
\lambda_j(\Delta'_{\eps, \delta})=\lambda_j(\Delta'_{\eps,
1})\approx\lambda_j(\Delta^\eps)\approx j^{\eps/2},
\end{equation}
where $\Delta$ is the usual Dirichlet Laplacian on $\{|w_1|<1,
|w_2|<1\}$ and the last estimate follows from the classical Weyl
formula.  Thus it follows from \eqref{v-1-estimate} that
\[
\|(\Delta'_{\eps, \delta})^{1/2}\chi_\delta
(I-P)\widetilde{N}_{\eps, \delta} u\|^2=\|\chi_\delta
v_1\|^{'2}_{\eps, \delta}\lesssim \|v_1\|^{'2}_{\eps,
\delta}\lesssim \|\widetilde{N}^{1/2}_{\eps, \delta} u\|^2_{\eps,
\delta}=\|u\|^2.
\]
Therefore, by \eqref{minmax-2} and \eqref{weyl}, we have
\begin{equation}\label{eigen-estimate-2}
\lambda_j(\chi_\delta(I-P)\widetilde{N}_{\eps, \delta})\lesssim
\lambda_j((\Delta'_{\eps, \delta})^{-1/2})\approx (1+j)^{-\eps/4}.
\end{equation}

We now study the eigenvalues of $\chi_\delta P\widetilde{N}_{\eps,
\delta}$.  Let $\scriptr\colon L^2(\C\times (0, \infty))\to H$ be
defined by
\[
\scriptr \phi(w_1, w_2)=\frac{1}{\sqrt{\pi}}\int_0^\infty
e^{w_2\nu}\phi (w_1, \nu)\sqrt{\nu}\, d\nu,
\]
and let $\scriptr^*\colon L^2(\C^2_-)\to L^2(\C\times (0,
\infty))$ be defined by
\[
\scriptr^*u(w_1, \nu)=\sqrt{\frac{\nu}{\pi}}\int_{-\infty}^0
(\scriptf_t u)(w_1, s, \nu) e^{s\nu}\, ds,
\]
where, as usual, $\scriptf_t$ is the Fourier transform in the
$t$-variable.  It is easy to see that $\scriptr$ is isometric and
onto.  Furthermore, $\scriptr^*\scriptr=I$ and
$\scriptr\scriptr^*=P$.

For any $\lambda>1$, let $\scripte_\lambda\colon L^2(\C\times (0,
\infty))\to L^2(\C\times (0, \infty))$ be defined by
\[
\scripte_\lambda \phi (w_1,
\nu)=(\scriptf^{-1}_{w_1}\chi_{\{1+|\zeta_1|^2+\nu^2<\lambda^2\}}
\scriptf_{w_1}\phi)(w_1, \nu),
\]
where $\chi_A$ is the characteristic function for the set $A$ as
before and $\scriptf_{w_1}$ is the Fourier transform in the $x$
and $y$ variables. (Recall that $w_1=x+iy$, $w_2=s+it$, and their
dual variables are $\zeta_1=\xi+i\eta$ and $\zeta_2=\mu+i\nu$.)
Let $\scriptm_\lambda=R\scripte_\lambda R^*\colon L^2(\C^2_-)\to
H$. Then $\scriptm_\lambda$ is an orthogonal project into $H$.  A
straightforward calculation yields that the kernel $M_\lambda$ of
$\scriptm_\lambda$ is given by
\[
M_\lambda(w, w')=\frac{1}{4\pi^3}\int_0^\infty\int_{\C} \nu
e^{(w_2+\ov{w'}_2)\nu
+i((x-x')\xi+(y-y')\eta)}\chi_{\{1+|\zeta_1|^2+\nu^2<\lambda^2\}}\,
d\xi d\eta d\nu.
\]
Thus,
\[
M_\lambda (w, w)=\frac{1}{4\pi^2}\int_0^{\sqrt{\lambda^2-1}}\nu
e^{2s\nu}(\lambda^2-1-\nu^2)\, d\nu .
\]
The square of the Hilbert-Schmidt norm of $\chi_\delta
\scriptm_\lambda$ equals
\[
\begin{aligned}
\int_{\C^4}|\chi_\delta(w)M_\lambda (w, w')|^2\, dV(w,w')
&=\int_{\C^2} |\chi_\delta|^2 M_\lambda (w, w)\, dV(w)\\
&=\frac{1}{4\pi^2\delta}\int_{\C^2} |\chi(w_1, w_2)|^2
M_\lambda(w_1, \delta^{-1/2}w_2)\, dV(w) \\
&\lesssim \delta^{-1}\int_{-\infty}^0\, ds\,
\int_0^{\sqrt{\lambda^2-1}}\nu
e^{2\delta^{-1/2}s\nu}(\lambda^2-1-\nu^2)\, d\nu \\
&\lesssim
\delta^{-1/2}\int_0^{\sqrt{\lambda^2-1}}(\lambda^2-1-\nu^2)\, d\nu
\lesssim\delta^{-1/2}\lambda^3.
\end{aligned}
\]
Thus, on the one hand, we have
\begin{equation}\label{eigen-estimate-3}
\lambda_j(\chi_\delta\scriptm_\lambda)\lesssim
(\delta^{-1/2}\lambda^3/j)^{1/2}.
\end{equation}
On the other hand, we have
\begin{align}
\|v_2-\scriptm_\lambda v_2\|^2 &=\|\scriptr\scriptr^*
v_2-\scriptr\scripte_\lambda\scriptr^*
v_2\|^2=\|(I-\scripte_\lambda)\scriptr^* v_2\|^2 \notag\\
&=\|(1-\chi_{\{1+|\zeta_1|^2+\nu^2<\lambda^2\}})\scriptf_{w_1}
\scriptr^* v_2\|^2 \notag \\
&\lesssim \int_0^\infty \, d\nu \int_{\C}
\chi_{\{1+|\zeta_1|^2+\nu^2\ge\lambda^2\}}\,d\xi d\eta
\left|\sqrt{\nu}\int_{-\infty}^0 (\scriptf_{w_1}\scriptf_t v_2)
(\zeta_1, \nu, s) e^{s\nu}\, ds\right|^2 \notag\\
&\lesssim\int_0^\infty \, d\nu \int_{\C}
\chi_{\{1+|\zeta_1|^2+\nu^2\ge\lambda^2\}}\,d\xi d\eta
\int_{-\infty}^0 |\scriptf_{w_1}\scriptf_t v_2|^2 \, ds
\notag\\
&\lesssim \lambda^{-2\eps}\vvv v_2 \vvv^2_{\eps}\lesssim
\lambda^{-2\eps}\|\widetilde{N}^{1/2}_{\eps, \delta} u\|^2_{\eps,
\delta}=\lambda^{-2\eps}\|u\|^2. \label{eigen-estimate-4}
\end{align}
(Here we have used Lemma~\ref{Bergman} in the last estimate.) From
\eqref{eigen-estimate-3}, \eqref{eigen-estimate-4}, and
\eqref{minmax-0}, we obtain
\begin{align}
 \lambda_j(\chi_\delta P\widetilde{N}^{1/2}_{\eps,
 \delta})&\lesssim (\delta^{-1/2}\lambda^3
 j^{-1})^{1/2}+\lambda^{-\eps} \notag \\
&\lesssim (j\delta^{1/2})^{-1/8}+(j\delta^{1/2})^{-\eps/4}\lesssim
(j\delta^{1/2})^{-\eps/4}, \label{eigen-esitmate-5}
\end{align}
by taking $\lambda=(j\delta^{1/2})^{1/4}$. This estimate, combined
with \eqref{eigen-estimate-2}, then gives us the desired estimate.
\end{proof}

\section{Estimates on the spectral kernel}\label{lower}

Let $E(\lambda)$ be the spectral resolution of the $\dbar$-Neumann
Laplacian $\square$ on $L^2_{(0, 1)}(\Omega)$ and let $e(\lambda; z,
\zeta)$ be its kernel.  By the classical elliptic theory we know
that
\begin{equation}\label{garding}
\lim_{\lambda\to\infty}\lambda^{-2}\trace e(\lambda; z,
z)=(2\pi)^{-2},
\end{equation}
where the limit is uniform on any compact subset of $\Omega$ ({\it
e.g.}, \cite{Garding53}). In fact, $\trace e(\lambda; z,
z)\lesssim \lambda^2$ for all $z\in\Omega$ with $d(z)\ge
\lambda^{-1/2}$ ({\it e.g.}, \cite{Metivier81}).

\begin{lemma}\label{kernel-estimate-1} Let $\tau=1/\sqrt{\lambda}$.
For sufficiently large $\lambda>0$,
\begin{equation}\label{kernel-estimate-2}
\trace e(\lambda; z, z)\lesssim \lambda(\delta(\pi(z),
\tau))^{-1},
\end{equation}
for all $z\in \Omega$ with $d(z)\gtrsim (\delta(\pi(z),
\tau))^{1/2}$.
\end{lemma}
\begin{proof} In light of the above remarks, it suffices to prove
\eqref{kernel-estimate-2} when $(\delta(\pi(z),
\tau))^{1/2}\lesssim d(z)\lesssim \lambda^{-1/2}$.  We will use a
global rescaling scheme which is slightly different from the local
rescaling scheme introduced in Section~\ref{rescaled}.

Let $z'\in b\Omega$ and let $\Omega'=\Phi'(\Omega)$ as in
Section~\ref{polydisc} where $\Phi'$ is given by
\eqref{coordinate-1}. Let $\delta=\delta(z', \tau)$.  For any
$\sigma$ such that $\sqrt{\delta}\lesssim \sigma\lesssim \tau$,
let $(\xi'_1, \xi'_2)=D'_{\tau, \sigma}(\xi_1, \xi_2)=(\xi_1/\tau,
\xi_2/\sigma)$. Let $\Omega'_{\tau, \sigma}=D'_{\tau,
\sigma}(\Omega)$. Let
\begin{equation}\label{l-prime}
L'_1=\frac{1}{|\partial r|}\left(\frac{\partial r}{\partial
z_2}\frac{\partial}{\partial z_1}-\frac{\partial r}{\partial
z_1}\frac{\partial}{\partial z_2}\right),  \qquad
L'_2=\frac{1}{|\partial r|}\left(\frac{\partial r}{\partial \bar
z_1}\frac{\partial}{\partial z_1}+\frac{\partial r}{\partial\bar
z_2}\frac{\partial}{\partial z_2}\right).
\end{equation}
We extend the vector fields $L'_1$ and $L'_2$ to form an
orthonormal basis for $T^{1, 0}(\C^2)$ over $\ov{\Omega}$. Let
$\omega'_1$ and $\omega'_2$ be the dual basis. Denote by
$L'^{\tau, \sigma}_1$, $L'^{\tau, \sigma}_2$, $\omega'^{\tau,
\sigma}_1$, and $\omega'^{\tau, \sigma}_1$ the vector fields and
forms defined as above but with $z$ replaced by $\xi'$ and $r$
replaced by $\rho_{\tau, \sigma}(\xi')=(1/\sigma)\rho(\tau \xi'_1,
\sigma\xi'_2)$ where $\rho$ is given by \eqref{defining-0}. Let
$\scripth_{\tau, \sigma}\colon L^2_{0, 1}(\Omega'_{\tau,
\sigma})\to L^2_{0, 1}(\Omega)$ be the unitary transformation
defined by
\[
\scripth_{\tau, \sigma}(v_1\ov{\omega}'^{\tau,
\sigma}_1+v_2\ov{\omega}'^{\tau, \sigma}_2)=
(\tau\sigma)^{-1}J\Phi'(z)(v_1(D'_{\tau,
\sigma}\circ\Phi')\ov{\omega}'_1+ v_2(D'_{\tau, \sigma}\circ
\Phi')\ov{\omega}'_2),
\]
where as before $J\Phi'$ is the Jacobian determinant of $\Phi'$.
Let
\[
Q^*_{\tau, \sigma}(u, v)=\tau^2 Q(\scripth_{\tau, \sigma}u,
\scripth_{\tau, \sigma}v)
\]
with $\dom(Q^*_{\tau, \sigma})=\dom(Q'_{\tau, \sigma})$ where $Q$ is
the sesquilinear form associated with the $\dbar$-Neumann Laplacian
$\square$ on $\Omega$ as before, and $Q'_{\tau, \sigma}$ is
associated with the $\dbar$-Neumann Laplacian $\square'_{\tau,
\sigma}$ on $\Omega'_{\tau, \sigma}$. Let $\square^*_{\tau, \sigma}$
be the self-adjoint operator defined by $Q^*_{\tau, \sigma}$ and let
$e^*_{\tau, \sigma}(\lambda; \xi', \eta')$ be the kernel of the
spectral resolution of $\square^*_{\tau, \sigma}$. Then
\begin{equation}\label{transformation-1}
e(\lambda; z, \zeta)=(\tau\sigma)^{-2} J\Phi' (z)
\ov{J\Phi'(\zeta)} e^*_{\tau, \sigma}(\tau^2\lambda; D'_{\tau,
\sigma}\circ\Phi'(z), D'_{\tau, \sigma}\circ\Phi'(\zeta)).
\end{equation}
Let $P'=\{\xi'\in \C^2 \mid |\xi'_1|<c, |\xi'_2+1/2|<c\}$. Then
for sufficiently small $c>0$, $P'$ is relatively compact subset of
$\Omega'_{\tau, \sigma}$. Furthermore, if $u$ is supported in
$P'$, then
\begin{equation}\label{transformation-2}
Q^*_{\tau, \sigma}(u, u)\gtrsim
\tau^2\|(\tau\sigma)^{-1}\nabla_{\xi} u(\tau^{-1}\xi'_1,
\sigma^{-1}\xi'_2)\|^2_{\Omega'}\gtrsim \|\nabla_{\xi'}
u\|^2_{\Omega'_\tau},
\end{equation}
where the rescaling of $u$ is component-wise and last estimate
follows from $\sigma\lesssim \tau$. In \eqref{transformation-1}, we
take $z=(\Phi')^{-1}(0, -\sigma/2)$. Then $d(z)\approx \sigma$. From
\eqref{transformation-1}, we have
\begin{equation}\label{transformation-3}
\trace e(\lambda; \ z, z)\lesssim (\tau\sigma)^{-2} \trace
e^*_{\tau, \sigma}(1; (0, -1/2), (0, -1/2)).
\end{equation}
By \eqref{transformation-2} and the Sobolev lemma,  we have for
any $k>2$ and $\xi'\in P'$
\begin{align*}
\|\trace e^*_{\tau, \sigma}(\tau^2\lambda; \cdot,
\xi')\|_{L^\infty(P')}&\lesssim \|\square^{*k}_{\tau, \sigma} \trace
e^*_{\tau, \sigma}(\tau^2\lambda; \cdot, \xi')\|_{L^2(\Omega'_{\tau,
\sigma})}+\|\trace e^*_{\tau, \sigma}
(\tau^2\lambda; \cdot, \xi')\|_{L^2(\Omega'_{\tau, \sigma})}\\
&\lesssim (1+(\tau^2\lambda)^k)(\trace e^*_{\tau,
\sigma}(\tau^2\lambda; \xi', \xi'))^{1/2}.
\end{align*}
Therefore $\trace e^*_{\tau, \sigma}(1; (0, -1/2), (0,
-1/2))\lesssim 1$. As $z'$ varies on $b\Omega$ and $\sigma$ varies
between $\sqrt{\delta}\lesssim \sigma\lesssim\tau$, we obtain
\eqref{kernel-estimate-2} from \eqref{transformation-3} for
$z\in\Omega$ with $\delta^{1/2}\lesssim d(z)\lesssim \tau$. This
concludes the proof of the lemma.
\end{proof}

We now estimate the spectral kernel $e(\lambda; z, z)$ for $z$
closed to the boundary using the rescaling scheme in
Section~\ref{rescaled}.  We will keep the notations of
Section~\ref{rescaled}. Let $\kappa$ be a cut-off function
compactly supported on $U_{z'}$ and identically 1 on a
neighborhood of $z'$ of uniform size.  Let
\[
E_\tau(\lambda)=\scriptg^{-1}_{\tau}\kappa
E(\lambda/\tau^2)\kappa\scriptg_{\tau}\colon
(L^2(\widetilde{\Omega}_{\tau}))^2\to
(L^2(\widetilde{\Omega}_{\tau}))^2.
\]
Then the kernel of $E_\tau(\lambda)$ is given by
\begin{equation}\label{spectral-1}
e_{\tau}(\lambda; w, w')=e(\lambda/\tau^2; \Psi^{-1}_{\tau} (w),
\Psi^{-1}_{\tau}(w'))\kappa(w)\kappa(w')|\det
d\Psi^{-1}_{\tau}(w)|^{\frac{1}{2}}|\det
d\Psi^{-1}_{\tau}(w')|^{\frac{1}{2}},
\end{equation}
where $\Psi_\tau=\widetilde\Psi_\tau\circ\widehat\Psi$.

\begin{lemma}\label{spectral-2}  Let
$\tau=1/\sqrt{\lambda}$.  Then for any $z'\in b\Omega$ and
sufficiently large $\lambda>0$,
\begin{equation}\label{spectral-3}
\int_{R_{\tau}(z')\cap\Omega} \trace e(\lambda; z, z)
dV(z)\lesssim (\delta(z', \tau))^{-1/2}.
\end{equation}
\end{lemma}

\begin{proof} In light of \eqref{spectral-1}, it suffices to prove
\begin{equation}\label{spectral-4}
\int_{P_{\tau}(z')\cap\widetilde{\Omega}_{\tau}} \trace
e_{\tau}(1; w, w)\, dV(w) \lesssim (\delta(z', \tau))^{-1/2}.
\end{equation}
Let $\square_{\tau}\colon (L^2(\widetilde{\Omega}_{\tau}\cap
P_\tau(z')))^2\to (L^2(\widetilde{\Omega}_{\tau}\cap P_\tau(z')))^2$
be the operator associated with the sesquilinear form $Q_{\tau}$
that is given by \eqref{Q-tau-def} and has domain
$\{\scriptg^{-1}_{\tau}(u) \mid u\in \dom(Q), \supp u\subset
\ov{\Omega}\cap \Psi^{-1}_\tau(P_\tau(z'))\}$. Note that
$\square_{\tau}=\tau^2 \scriptg_{\tau}^{-1}\square \scriptg_{\tau}$.
Also by Lemma~\ref{phi-tilde}, $R_{C^{-1}\tau}(z')\subset
\Psi^{-1}_\tau(P_\tau(z'))\subset R_{C\tau}(z')$ for some
sufficiently large constant $C>0$.  By Lemma~\ref{Q-estimate} and
using the fact that $\delta=\delta(z', \tau)\lesssim \tau^2$, we
have
\[
Q_\tau(u, u)\gtrsim \|u\|^2_{\eps, C\delta},
\]
for $u\in \dom(\square^{1/2}_\tau)$, where $C>0$ is any constant.
Therefore,
\[
\|u\|^2=Q_\tau(N^{1/2}_\tau u, N^{1/2}_\tau u)\gtrsim
\|\widetilde{\square}^{1/2}_{\eps, C\delta} N^{1/2}_\tau u\|^2.
\]
Choose $C>2$.  Then $N^{1/2}_\tau=\chi_{C\delta}N^{1/2}_\tau
=\chi_{C\delta}\widetilde{N}^{1/2}_{\eps,
C\delta}\widetilde{\square}^{1/2}_{\eps, C\delta}N^{1/2}_\tau$. It
follows from \eqref{minmax-2} and Lemma~\ref{eigen-estimate-1}
that
\begin{equation}\label{n-tau}
\lambda_j(N_{\tau}^{1/2}) \lesssim
\lambda_j(\chi_{C\delta}\widetilde{N}_{\eps, \delta}^{1/2})\lesssim
(j\delta^{1/2})^{-\eps/4}.
\end{equation}

Let $K$ be any positive integer such that $K>4/\eps$.  Let
$\chi^{(k)}$, $k=0, 1, \ldots, K$, be a family of cut-off
functions supported in $\{|w_1|<1,  |w_2|<1\}$ such that
$\chi^{(0)}=\chi$ and $\chi^{(k+1)}=1$ on $\supp \chi^{(k)}$. Let
\[
E^{(l)}_\tau(\lambda)=\scriptg^{-1}_\tau \kappa (\tau^2\square)^l
E(\lambda\tau^{-2})\kappa\scriptg_\tau\colon
(L^2(\widetilde{\Omega}_\tau))^2\to
(L^2(\widetilde{\Omega}_\tau))^2.
\]
Note that $E^{(0)}_\tau(\lambda)=E_\tau(\lambda)$ and
\begin{equation}\label{e-ell}
\|E^{(l)}_\tau (1) u\|\lesssim \|u\|.
\end{equation}
It is easy to check the following commutating identity:
\begin{equation}\label{Q-commute}
Q(\theta u, \theta u)=\Re (\theta \square u, \theta u)+(1/2)(u,
[\theta, A] u)
\end{equation}
where $\theta$ is any smooth function and $A=[\dbar^*,
\theta]\dbar+\dbar [\dbar^*, \theta]+\dbar^*[\dbar,
\theta]+[\dbar, \theta]\dbar^*$.  Note that $[\theta, A]$ is of
zero order.  Using the above identity and the Schwarz inequality,
we obtain that for any $u\in (L^2(\widetilde{\Omega}_{\tau}))^2$,
\begin{align*}
\|\square_{\tau}^{\frac{1}{2}} \chi^{(k)}_\delta E^{(l)}_\tau (1)
u\|^2 &=Q_{\tau}(\chi^{(k)}_\delta
E^{(l)}_\tau(1)u, \; \chi^{(k)}_\delta E^{(l)}_\tau(1)u)\\
&=\tau^2 Q(\chi^{(k)}_\delta(\Psi_{\tau}) (\tau^2\square)^l
E(\tau^{-2})\kappa\scriptg_{\tau} u,\;
\chi^{(k)}_\delta(\Psi_{\tau}) (\tau^2\square)^l
E(\tau^{-2})\kappa\scriptg_{\tau} u)\\
&\lesssim \|\chi^{(k)}_\delta(\Psi_{\tau}) (\tau^2\square)^{l+1}
E(\tau^{-2})\kappa\scriptg_{\tau} u\|^2+
\|\chi^{(k+1)}_\delta(\Psi_{\tau}) (\tau^2\square)^l
E(\tau^{-2})\kappa\scriptg_{\tau} u\|^2 \\
&=\|\chi^{(k)}_\delta E^{(l+1)}_\tau (1)
u\|^2+\|\chi^{(k+1)}_\delta E^{(l)}_\tau (1) u\|^2.
\end{align*}
By \eqref{minmax-1} and \eqref{minmax-2}, we have
\begin{equation}\label{lambda-j}
\begin{aligned}
\lambda_{3j+1}(\chi^{(k)}_\delta E^{(l)}_\tau (1))&\le
\lambda_{j+1} (N_\tau^{1/2})\lambda_{2j+1}(\square^{1/2}_\tau
\chi^{(k)}_\delta E^{(l)}_\tau (1))\\
&\le \lambda_{j+1}
(N_\tau^{1/2})\left(\lambda_{j+1}(\chi^{(k)}_\delta E^{(l+1)}_\tau
(1))
+\lambda_{j+1}(\chi^{(k+1)}_\delta E^{(l)}_\tau(1))\right).\\
\end{aligned}
\end{equation}
Using \eqref{n-tau}, \eqref{e-ell}, and\eqref{lambda-j}, we then
obtain by an inductive argument on $K-(k+l)$ that
\begin{equation}\label{lambda-e-ell}
\lambda_j(\chi^{(k)}_\delta E^{(l)}_\tau (1))\lesssim
(j\delta^{1/2})^{-(K-(k+l))\eps/4}
\end{equation}
for any pair of non-negative integers $k, l$ such that $0\le
k+l\le K$. In particular,
\[
\lambda_j(\chi_\delta E_\tau(1))\lesssim (j\delta^{1/2})^{-K\eps/4}.
\]
Since $E_\tau(1)$ is a contraction, we also have that
$\lambda_j(\chi_\delta E_\tau(1))\lesssim 1$.  The trace norm of
$\chi_\delta E_\tau (1)$ is then given by
\[
\sum_{j\le \delta^{-1/2}} \lambda_j(\chi_\delta
E_\tau(1))+\sum_{j>\delta^{-1/2}} \lambda_j(\chi_\delta E_\tau(1))
\lesssim \delta^{-1/2}+ \sum_{j>\delta^{-1/2}}
(j\delta^{-1/2})^{-K\eps/4} \lesssim \delta^{-1/2}.
\]
Inequality \eqref{spectral-4} is now an easy consequence of the
above estimate.
\end{proof}

\begin{lemma}\label{spectral-5} For sufficiently large
$\lambda>0$,
\begin{equation}\label{spectral-6}
\int_{A_\tau} \trace e(\lambda; z, z)\, dV(z)\lesssim
\tau^{-2}\int_{b\Omega} (\delta(z', \tau))^{-1}\, dS(z')\lesssim
\lambda^{1+m}.
\end{equation}
\end{lemma}

\begin{proof} By Lemma~\ref{spectral-2}, we have
\[
(\delta(z', \tau))^{-1/2}\int_{\Omega}\chi_{\Omega\cap
R_\tau(z')}(z) \trace e(\lambda; z, z)\, dV(z)\lesssim (\delta(z',
\tau))^{-1}.
\]
Integrating both sides with respect to $z'\in b\Omega$ and using
the Fubini-Tonelli Theorem, we have
\[
\int_{\Omega} \trace e(\lambda; z, z)\, dV(z) \int_{z'\in b\Omega}
\chi_{\Omega\cap R_\tau(z')}(z) (\delta(z', \tau))^{-1/2}\,
dS(z')\lesssim \int_{b\Omega} (\delta(z', \tau))^{-1}\, dS(z').
\]
By Lemma~\ref{delta-int} and the fact that $\delta(z',
\tau)\gtrsim \tau^{2m}$, we then obtain
\[
\int_{A_{C\tau}} \trace e(\lambda; z, z)\, dV(z)\lesssim
\tau^{-2}\int_{b\Omega} (\delta(z', \tau))^{-1}\, dS(z')\lesssim
\lambda^{1+m}.
\]
Lemma~\ref{spectral-5} then follows from a rescaling of $\tau$ in
the above arguments.
\end{proof}

We are now in position to prove the following variation of
Theorem~\ref{maintheorem2a}.

\begin{proposition}\label{maintheorem2}
Let $\Omega$ be a smooth bounded pseudoconvex domain of finite
type $2m$. Then
\begin{equation}\label{lower-estimate}
\limsup_{\lambda\to\infty}
\frac{\scriptn(\lambda)}{\lambda^{m+1}}\lesssim
\limsup_{\lambda\to\infty} \frac{1}{\lambda^m}\int_{b\Omega}
(\delta(z', 1/\sqrt{\lambda}))^{-1}\, dS(z')\lesssim 1.
\end{equation}
\end{proposition}

\begin{proof} Note that
\[
\begin{aligned}
\frac{\scriptn(\lambda)}{\lambda^{1+m}}&=\lambda^{-1-m}\int_{\Omega}\trace e(\lambda; z, z)\, dV(z)\\
          &=\lambda^{-1-m}\int_{A_\tau} \trace e(\lambda; z, z)\, dV(z)+
          \lambda^{-1-m}\int_{\Omega\setminus A_\tau}
          \trace e(\lambda; z, z)\, dV(z).
\end{aligned}
\]
By Lemma~\ref{spectral-5}, the first term in the last expression
is bounded above by
\[
\lambda^{-m}\int_{b\Omega} (\delta(z', 1/\sqrt{\lambda}))^{-1}\,
dV(z).
\]
By \eqref{garding}, Lemma~\ref{kernel-estimate-2}, and the
Lebesgue dominated convergence theorem, we have
\[
\lim_{\lambda\to\infty} \lambda^{-1-m}\int_{\Omega\setminus
A_\tau} \trace e(\lambda; z, z)\, dV(z) =
\begin{cases}
(2\pi)^{-2}\text{vol}(\Omega), &\qquad \text{if} \ m=1\\
0, &\qquad \text{if} \  m>1.
\end{cases}
\]
We then conclude the proof of the proposition by noting that
$\lambda^{-m}\lesssim \delta(z', 1/\sqrt{\lambda})\lesssim
\lambda^{-1}$.
\end{proof}

\noindent{\bf Remarks.}  It follows from the above proof that
$\limsup_{\lambda\to\infty} \scriptn(\lambda)/\lambda^{m+1}=0$
when $m>1$.

\section{Hearing a finite type property}\label{proof1}

We prove Theorem~\ref{maintheorem3} in this section. The following
lemma is well-known (see \cite{FornaessSibony89} for the two
dimensional case and \cite{Yu95} for the general case). It can be
proved along the lines of the arguments in Section~\ref{polydisc}.
Throughout this section, we will use $z'$ to denote the first
$(n-1)$-tuple of $z\in\C^n$.

\begin{lemma}\label{multitype}
Let $\Omega$ be a smooth bounded pseudoconvex domain in $\C^n$.
Assume that the $D_{n-1}$-type of $b\Omega$ at $z^0$ is $\ge 2m$.
Then there exists a neighborhood $U$ of $z^0$ and a biholomorphic
map $w=\Psi(z)$ from $U$ into $\C^n$ such that $\Psi(z^0)=0$ and
\[
\Psi(\Omega\cap U)=\{w\in\C^n \mid |w'|<1, |\Im w_n|<1, \rho(w)=\Re
w_n+h(w', \Im w_n)<0\},
\]
where $h(w', \Im w_n)=f(w')+(\Im w_n)\cdot g(w') + \sigma(w', \Im
w_n)$ with $|f(w')|\lesssim |w'|^{2m}$, $|g(w')|\lesssim
|w'|^{m+1}$, and $|\sigma(w', \Im w_n)|\lesssim (\Im w_n)^2$.
\end{lemma}

We now prove Theorem~\ref{maintheorem3}.  Assume that the
$D_{n-1}$-type of $b\Omega$ at $z^0$ is $\ge 2m$.  We apply
Lemma~\ref{multitype} and keep its notations.

Write $w_n=s+it$.  Let $b(t)$ be the cut-off function constructed in
the paragraph preceding Lemma~\ref{wavelett}. We first extend $b(t)$
to the whole complex plane as follows.  Let $\chi$ be any smooth
cut-off function supported on $(-2,\ 2)$ and identically 1 on $(-1,
1)$ and let
\[
B(w_n)=(b(t)-ib'(t)s-b''(t)s^2/2)\chi(s/(1+|t|^2)).
\]
Then $B(0,t)=b(t)$ and $|\partial B(w_n)/\partial\bar w_n|\lesssim
|s|^2$.

Let
\[
L_j=\frac{\partial\rho}{\partial w_n}\frac{\partial}{\partial w_j}
-\frac{\partial\rho}{\partial w_j}\frac{\partial}{\partial w_n},\ \
1\le j\le n-1, \text{ and } L_n=\sum_{j=1}^n
\frac{\partial\rho}{\partial\bar w_j}\frac{\partial}{\partial w_j}.
\]
Let
\[
\tilde L_j=(\Psi^{-1})_*(L_j), \ \ 1\le j\le n-1, \text{ and }\tilde
L_n=\sum_{j=1}^n\frac{\partial\tilde\rho}{\partial\bar
z_j}\frac{\partial}{\partial z_j},
\]
where $\tilde\rho(z)=\rho\circ\Psi (z)$. Then $\tilde L_j, 1\le j\le
n$, form a basis for $T^{1, 0}(\C^n)$ in a neighborhood of $z^0$.
Replacing $\tilde L_j$ by the product of $\tilde L_j$ with an
appropriate cut-off function, we may assume that $\tilde L_j$ is
supported in $U$. Let $\hat L_j$, $1\le j\le n$, be the vector
fields obtained after performing the Gram-Schmidt process on $\tilde
L_j$. Let $\hat \omega_j$, $1\le j\le n$, be the dual basis of $\hat
L_j$. Let $a(w')$ be a smooth function identically 1 near the origin
and compactly supported in the unit ball in $\C^{n-1}$.  For any
positive integers $j$ and for any positive integer $k$ such that
$2^{mj-1}/j\le k\le 2^{mj}/j$, let
\[
f_{j, k}(w)=k8^{(m+n-1)j} a(8^jw')B(8^{mj}w_n)e^{2\pi k^28^{mj}w_n}.
\]
Let $g_{j, k}(z)=f_{j, k}(\Psi(z))\cdot (J\Psi(z))$.  Let
\[
u_{j, k}=g_{j, k}(z)
\bar{\hat{\omega}}_1\wedge\cdots\wedge\bar{\hat{\omega}}_{n-1}.
\]
Then for any sufficiently large $j$, $u_{j, k}$ is a compactly
supported smooth $(0, n-1)$-form in $\dom(Q_{n-1})$. Moreover,
\begin{align*}
\|u_{j, k}\|^2_\Omega &= \int_{\C^{n-1}} dV(w')\int_\R dt
\int^{-h(w', t)}_{-\infty} |f_{j, k}(w)|^2 \,ds \\
&=k^2\int_{\C^{n-1}} |a(\tilde w')|^2 dV(\tilde w')\int_\R \,d\tilde
t\int^{-8^{mj}h(8^{-j}\tilde{w}', 8^{-mj}\tilde{t})}_{-\infty}
|B(\tilde s, \tilde t)|^2 e^{4\pi k^2 \tilde s}\, d\tilde s.
\end{align*}
(After the substitution $\tilde w'=8^j w',\ \tilde w_n=8^{mj}w_n$.)
Since $|8^{mj}h(8^{-j}\tilde w', 8^{-mj}\tilde t)|\lesssim 8^{-mj}$
and $k^2 8^{-mj}\le j^{-2}2^{-mj}$, we have
\begin{align*}
\|u_{j, k}\|^2_\Omega &\lesssim k^2\int_{\C^{n-1}} |a(\tilde w')|^2
dV(\tilde w')\int_\R \,d\tilde t\int^{C 8^{-mj}}_{-\infty} |B(\tilde
s, \tilde t)|^2 e^{4\pi k^2
\tilde s}\, d\tilde s\\
&\lesssim \int_{\C^{n-1}}|a(\tilde w')|^2 dV(\tilde w')\int_{|\tilde
t|<1} e^{Ck^2 8^{-mj}} d\tilde t \lesssim 1.
\end{align*}
Similarly, $\|u_{j, k}\|^2\gtrsim 1$.  Therefore, $\|u_{j,
k}\|^2_\Omega \approx 1$. Furthermore, after a substitution as
above, we have that for any $k$, $l$ such that $2^{mj-1}/j\le k, l
\le 2^{mj}/j$, $\langle u_{j, k}, u_{j, l}\rangle $ equals
\[
kl \int_{\C^n} |a(\tilde w')|^2 dV(\tilde w')\int_\R d\tilde t
\int^{-8^{mj}h(8^{-j}\tilde w', 8^{-mj}\tilde t)}_{-\infty}
|B(\tilde s, \tilde t)|^2 e^{2\pi((k^2+l^2)\tilde s+i(k^2-l^2)\tilde
t)}\, d\tilde s
\]
Let $A$ be the above expression  with the upper limit in the last
integral replaced by $0$ and let $B$ likewise be the above
expression with the lower limit of the last integral replaced by
$0$. Thus $\langle u_{j, k}, u_{j, l}\rangle=A+B$.  It is easy to
see that
\begin{align*}
|B| &\le kl \int_{\C^{n-1}} |a(\tilde w')|^2 dV(\tilde w')\int_\R
d\tilde t \int^{-8^{mj}h(8^{-j}\tilde w', 8^{-mj}\tilde t)}_0
|B(\tilde s, \tilde t)|^2
e^{2\pi(k^2+l^2)\tilde s}\, d\tilde s \\
&\lesssim \frac{kl}{k^2+l^2}\left(1-e^{C(k^2+l^2)8^{-mj}}\right)
\lesssim j^{-2}2^{-mj}.
\end{align*}
To estimate $|A|$, we first observe that by Lemma~\ref{wavelett},
for $k\not=l$,
\[
A=kl\int_{\C^{n-1}} |a(\tilde w')|^2 dV(\tilde w') \int_\R d\tilde
t\int^0_{-\infty} \left(|B(\tilde s, \tilde t)|^2-|B(0, \tilde
t)|^2\right)e^{2\pi((k^2+l^2)\tilde s+i(k^2-l^2)\tilde t)} \,d\tilde
s.
\]
Hence
\[
|A|\lesssim kl\int_{|\tilde w'|<1} |a(\tilde w')|^2 dV(\tilde w')
\int_{-1}^1\,d\tilde t\int^0_{-\infty} \tilde s
e^{2\pi(k^2+l^2)\tilde s}\, d\tilde s \lesssim kl/(k^2+l^2)^2.
\]
Therefore, for sufficiently large $j$ and for any $k$, $l$ such that
$2^{mj-1}/j\le k, l \le 2^{mj}/j$, $k\not= l$, we have,
\[
|\langle u_{j, k},\ u_{j, l}\rangle |\lesssim j^{-2}2^{-mj}.
\]
For any $k$ such that $2^{mj-1}/j\le k\le 2^{mj}/j$ and for any
$c_k\in\C$, we have
\begin{align*}
\|\sum_k c_k u_{j, k}\|^2 &=\sum_{k} |c_k|^2\|u_{j,
k}\|^2-\sum_{\substack{k, l\\ k\not=l}}
c_k\overline{c}_l \langle u_{j, k}, \ u_{j, l}\rangle \\
&\ge \sum_{k} |c_k|^2\|u_{j, k}\|^2-j^{-2}2^{-mj}\big|\sum_{k}
c_k\big|^2 \\
&\gtrsim (1-j^{-4})\sum_{k} |c_k|^2 \gtrsim \sum_{k} |c_k|^2,
\end{align*}
where the summations are taken over all integers between
$2^{mj-1}/j$ and $2^{mj}/j$.

Since each $\hat L_k$, $1\le k\le n-1$, is a linear combination of
$\tilde L_1$, \ldots, $\tilde L_k$, and $\hat L_n$ is just the
normalization of $\tilde L_n$, it follows that
\begin{align*}
Q_{n-1}(u_{j, k}, u_{j, k})& \lesssim \|g_{j,
k}\|^2_{\Omega}+\sum_{l=1}^{n-1} \big(\|\overline{\tilde L}_l g_{j,
k}\|^2_{\Omega}+\|\tilde L_l g_{j, k}\|^2_{\Omega}\big) +
\|\overline{\tilde L}_n g_{j,
k}\|^2_{\Omega}\\
&\lesssim \|f_{j, k}\|^2_{\Psi(\Omega\cap U)}+\sum_{l=1}^{n-1} \|L_l
f_{j, k}\|^2_{\Psi(\Omega\cap U)} + \|\dbar f_{j,
k}\|^2_{\Psi(\Omega\cap U)}.
\end{align*}
For $1\le l \le n-1$,  $\|L_l f_{j, k}\|^2_{\Psi(\Omega\cap U)}$ is
bounded above by
\begin{align*}
&2\big(\big\|\frac{\partial\rho}{\partial w_n}\frac{\partial f_{j,
k}}{\partial w_l}\big\|^2_{\Psi(\Omega\cap
U)}+\big\|\frac{\partial\rho}{\partial w_l}\frac{\partial
f_{j, k}}{\partial w_n}\big\|^2_{\Psi(\Omega\cap U)}\big)\\
&\quad \lesssim 8^{2j}+\int_{|\tilde w'|<1} |a(\tilde w')|^2
dV(\tilde w')\int_{-1}^1\, d\tilde t \int^{-8^{mj} h(8^{-j}\tilde
w', 8^{-mj}\tilde t)}_{-\infty} \Big( \big(8^{mj}|\nabla
B|+\\
&\qquad k^2 8^{mj}|B|\big)\big(|8^{-j}\tilde
w'|^{2m-1}+|8^{-mj}\tilde t|\cdot|8^{-j}\tilde w'|^{m-1}+
|8^{-mj}\tilde t|^2\big)\Big)^2  k^2 e^{4\pi k^2 \tilde s}\, d\tilde s\\
&\quad\lesssim 8^{2j}.
\end{align*}
Furthermore, $\|\dbar f_{j, k}\|^2_{\Psi(\Omega\cap U)}$ is
\begin{align*}
&\lesssim\sum_{l=1}^n \big\|\frac{\partial f_{j, k}}{\partial \bar
w_l}\big\|^2_{\Psi(\Omega\cap U)}\\
&\lesssim 8^{2j} +8^{2mj}k^2\int_{|\tilde w'|<1}|a(\tilde
w')|^2\,dV(\tilde w')\int_{-1}^1\,d\tilde t \int^{-8^{mj}
h(8^{-j}\tilde w', 8^{-mj}\tilde t)}_{-\infty} \big|\frac{\partial
B}{\partial \bar w_n}\big|^2 e^{4\pi k^2
\tilde s}\,d\tilde s\\
&\lesssim 8^{2j}+8^{2mj}k^2\int_{|\tilde w'|<1}|a(\tilde
w')|^2\,dV(\tilde w')\int_{-1}^1\,dt \int^{-8^{mj} h(8^{-j}\tilde
w', 8^{-mj}\tilde t)}_{-\infty} \tilde s^4 e^{4\pi k^2
\tilde s}\,d\tilde s\\
&\lesssim 8^{2j}+8^{2mj}k^{-8}\lesssim 8^{2j}.
\end{align*}
Therefore, we have
\[
Q_{n-1}(u_{j, k}, u_{j, k})\lesssim 8^{2j}.
\]
We now invoke the hypothesis of Theorem~\ref{maintheorem3}.  Since
$\scriptn_q(\lambda)$ has at most polynomial growth,
$\lambda_{j}(\square_q)\gtrsim j^\e$ for some $\e>0$. It follows
from Proposition~\ref{qeigenvalue} that
$\lambda_{j}(\square_{n-1})\gtrsim j^\e$. By Lemma~\ref{spectral},
for all sufficiently large $j$, there exists an integer $k_0\in
[2^{mj-1}/j, \ 2^{mj}/j]$ such that
\[
Q_{n-1}(u_{j, k_0}, u_{j, k_0})\gtrsim (2^{mj}/j)^\e.
\]
Therefore, $8^{2j} \gtrsim (2^{mj}/j)^\e$.  Hence $m\le 6/\e$. We
thus conclude the proof of Theorem~\ref{maintheorem3}.

\bibliography{survey}

\begin{thebibliography}{XXXXXX}
\bibitem[Ma97]{Matheos97}
P.~Matheos, \emph{A {H}artogs domain with no analytic discs in the
boundary for which the $\dbar$-{N}eumann problem is not compact},
UCLA Ph.D.\ thesis,  1997.
\end{thebibliography}
% Bibliography generated by BibTeX with amsplain
% bibliographystyle and mrabbrev.bib abbreviations.
%
%%%

%% bibliography generated by BiBTeX
\providecommand{\bysame}{\leavevmode\hbox to3em{\hrulefill}\thinspace}

\end{document}